\pdfoutput=1
\documentclass{article}

\usepackage{amsthm,amsmath,amssymb}
\usepackage{xcolor}
\usepackage{tikz}
\usepackage{hyperref}
\usepackage{mathrsfs}
\usepackage{comment}
\usepackage{lineno}

\newtheorem{thm}{Theorem}[section]

\newtheorem{lm}[thm]{Lemma}
\newtheorem{clm}[thm]{Claim}
\newtheorem{col}[thm]{Corollary}

\newtheorem{cons}[thm]{Construction}
\newtheorem{prop}[thm]{Proposition}
\newtheorem{conj}[thm]{Conjecture}
\def\eps{\varepsilon}
\def\I{\mathcal{I}}

\usepackage[hmargin=2.5cm,vmargin=3cm,bindingoffset=0.5cm]{geometry}
\title{Number of independent transversals in multipartite graphs}

\author{
Yantao Tang
\thanks{Department of Mathematics and Statistics, Georgia State University, Atlanta, GA 30303. Email: \texttt{ytang26@gsu.edu}.}
\and 
Yi Zhao
\thanks{Department of Mathematics and Statistics, Georgia State University, Atlanta, GA 30303. Email: \texttt{yzhao6@gsu.edu}.}
}
\date{}

\begin{document}
\maketitle

\begin{abstract}
An \emph{independent transversal} in a multipartite graph is an independent set that intersects each part in exactly one vertex. We show that for every even integer $r\ge 2$, there exist $c_r>0$ and $n_0$ such that 
every $r$-partite graph with parts of size $n\ge n_0$ and maximum degree at most $rn/(2r-2)-t$, where $t=o(n)$, 
contains at least $c_r t n^{r-1}$ independent transversals. This is best possible up to the value of $c_r$. Our result confirms a conjecture of Haxell and Szab\'o from 2006 and partially answers a question raised by Erd\H{o}s in 1972 and studied by Bollob\'as, Erd\H os and Szemer\'edi in 1975.
 
We also show that, given any integer $s\ge 2$ and even integer $r\ge 2$, there exist $c_{r,s}>0$ and $n_0$ such that every $r$-partite graph with parts of size $n\ge n_0$ and maximum degree at most 
$rn/(2r-2)- c_{r, s} n^{1-1/s}$ contains an independent set with exactly $s$ vertices in each part. This is best possible up to the value of $c_{r, s}$ if a widely believed conjecture for the Zarankiewicz number holds. Our result partially answers a question raised by Di Braccio and Illingworth recently.
\end{abstract}


\section{Introduction}
Let $G$ be an $r$-partite graph with parts $V_1,\cdots,V_r$. An \textit{independent transversal} (IT) of $G$ is an independent set containing exactly one vertex from each $V_i$. Independent transversals have found applications to SAT, linear arboricity, hypergraph matchings, list coloring, and many other areas (see, for example, \cite{MR4152571} and references therein). It is NP-complete to decide whether an $r$-partite graph has an IT, and therefore natural to seek sufficient conditions for the existence of an IT.

Denote by $\mathcal{G}_r(n)$ the class of all $r$-partite graphs with parts of size $n$. Let $\Delta_r(n)$ be the largest integer $C$ such that every graph $G\in \mathcal{G}_r(n)$ with maximum degree less than $C$ has an independent transversal. The study of $\Delta_r(n)$ can be traced back to 1972, when Erd\H{o}s \cite[Problem 2]{MR0329905} asked for the complementary function $\delta_r(n)= (r-1)n - \Delta_r(n)$, which is the smallest integer $c$ such that every graph $G\in \mathcal{G}_r(n)$ with minimum degree greater than $c$ contains a complete graph on $r$ vertices.\footnote{Erd\H{o}s initially conjectured that $\delta_r(n)=(r-2)n$ but Seymour found counterexamples for $r\ge 4$, see \cite{BES1975}.}  
Let $\Delta_r =\lim_{n\rightarrow \infty}\frac{\Delta_r(n)}{n}$ (it is easy to see the limit exists, see \cite{BES1975}). 
In 1975, Bollob\'as, Erd\H os and Szemer\'edi \cite{BES1975} 
showed $\frac{2}{r}\le \Delta_r\le \frac{1}{2}+\frac{1}{r-2}$ and conjectured $\lim_{r\to \infty} \Delta_r = \frac{1}{2}$. This was confirmed by Haxell \cite{H2001} in 2001. After the work of Graver(cf. \cite{BES1975}) on $r=3$ and Jin \cite{J1992} on $r=4, 5$, the value of $\Delta_r(n)$ was completely determined by Szab\'o and Tardos \cite{ST2006} (for even $r$) and Haxell and Szab\'o \cite{HS2006} (for odd $r$). 

\begin{thm}[\cite{HS2006, ST2006}]\label{thm:threshold}
For all integers $n\ge 1$ and $r\ge 2$, 
\[\Delta_r(n) = 
\begin{cases}
\lceil \frac{r-1}{2(r-2)}n \rceil & \text{if $r$ is odd} \\
\lceil \frac{r}{2(r-1)}n \rceil & \text{if $r$ is even}.
\end{cases}
\]
\end{thm}

Erd\H{o}s \cite[Problem 2, Page 353]{MR0329905} wrote: 
    \emph{Several further problems can be raised. First of all, how many complete subgraphs of $r$ vertices must $G$ contain? This is perhaps not quite simple even for $r = 3$.}
Here $G$ refers to any graph $G\in \mathcal{G}_r(n)$ with $\Delta(G)< \Delta_r(n)$. Formally, given positive integers $n, r, t$, denote by $f_r(n,t)$ the largest integer $f$ such that every graph $G\in \mathcal{G}_r(n)$ with maximum degree at most $\Delta_r(n)-t$ has at least $f$ independent transversals. 
It is easy to see that $f_2(n, t)= tn$ for all $t\le n=\Delta_2(n)$.\footnote{Every bipartite graph $G\in G_2(n)$ with $\Delta(G)\le n-t$ contains at least $tn$ ITs (nonadjacent crossing pairs) while an $(n-t)$-regular graph $G\in G_2(n)$ has exactly $tn$ ITs.}

Bollob\'as, Erd\H os and Szemer\'edi \cite{BES1975} answered Erd\H{o}s' question for $r=3$ by showing $f_3(n, 1)=4$ for all $n\ge 4$. Furthermore, they proved that $t^3\le f_3(n,t)\le 4t^3$ for all $t\le \frac{n}{5}$. Later Jin \cite{J1998} showed that $\frac{1}{324}n^3\le f_4(n,1)\le \frac{8}{9}(n+1)^3$ for $n\ge 9$. Bollob\'as, Erd\H os and Szemer\'edi \cite{BES1975} also gave the following supersaturation result for all $r$ when $t$ is a linear function of $n$.
\begin{thm}[\cite{BES1975}]\label{thm:BES}
	Given $r\ge 2$ and $\varepsilon>0$, there exist $\delta>0$ and integer $n_0$ such that $f_r(n,\varepsilon n)\ge \delta n^r$ for all $n\ge n_0$. 
\end{thm}
The authors of \cite{BES1975} further wrote: \emph{We do not obtain interesting results for $f_r(n,t)$ when $t=o(n)$ for $r\ge 4$ though we believe they exist.}
Possibly motivated by Jin's result \cite{J1998} on $f_4(n,1)$, Haxell and Szab\'o \cite{HS2006} conjectured that $f_r(n,1)=\Theta(n^{r-1})$ for \emph{all even} $r$ (see Section~\ref{sec:cr} for their comments on odd $r$). 

In this paper we prove that $f_r(n,t)=\Theta(t n^{r-1})$ for all $t\ge 1$ when $r$ is even thus confirming the conjecture of Haxell and Szab\'o and partially answering the question of Erd\H{o}s.
Since Theorem~\ref{thm:BES} covers the case when $t=\Omega(n)$, we assume $t=o(n)$ in the following theorem.
\begin{thm}\label{mainresult}
		For any even integer $r\ge 4$, there exist constants $\zeta$ and $n_0$ such that for any $n\ge n_0$ and $1\le t\le \zeta n$, we have
		\[
		\frac{1}{4}\left((t-1)(2r-2)+1\right)\left(\frac{r n}{2r-2}\right)^{r-1}\le f_r(n,t)\le r^2 t \left(\frac{rn}{2r-2}\right)^{r-1}.		
		\]
\end{thm}
We can replace $\frac{1}{4}$ by $\frac{1}{2} - o(1)$ in the lower bound and reduce $r^2$ in the upper bound by more careful calculations but choose the current form for simplicity. 
Nevertheless, Theorem~\ref{mainresult} shows that $\frac{1}{4}\left(\frac{r n}{2r-2}\right)^{r-1}\le f_r(n,1)\le r^2 \left(\frac{rn}{2r-2}\right)^{r-1}$ for $n\ge n_0$. In particular, $f_4(n,1)\ge \frac{2}{27}n^3$, better than the bound in \cite{J1998}. 

\medskip
Bollob\'as, Erd\H os and Szemer\'edi \cite{BES1975} also studied the minimum degree in a tripartite graph $G\in \mathcal{G}_3(n)$ that guarantees a copy of $K_3(2)$ (complete tripartite graph with 2 vertices in each part). They conjectured that $\delta(G)\ge n+ C n^{1/2}$ (some $C>0$) suffices. Here the first term is $\delta_3(n)=2n - \Delta_3(n)=n$ and the second term $n^{1/2}$ is necessary; indeed, there are many constructions of $G$ with $\delta(G)= n+ \Omega(n^{1/2})$ and without $K_3(2)$, see \cite{BZ2022,CHEN2027115475,MR5101285}. After two recent works \cite{BZ2022,CHEN2027115475}, this problem was settled very recently by Di Braccio and Illingworth \cite{MR5101285}, who showed that every $G\in \mathcal{G}_3(n)$ with $\delta(G)\ge n+ C n^{1- 1/s}$ contains a copy of $K_3(s)$. They also showed that $n^{1- 1/s}$  is best possible by assuming a well-known conjecture $z(n; s)=\Omega(n^{2- 1/s})$, where \emph{Zarankiewicz number} $z(n;s)$ is the largest number of edges in a bipartite graph $G\in \mathcal{G}_2(n)$ that does not contain a copy of $K_{s,s}$ (complete bipartite graph with $s$ vertices in each part). 

%I add a \red here to show the difference with the submitted version, but I didn't change any word here.
Di Braccio and Illingworth \cite[Problem 6.3]{MR5101285} asked the corresponding problem for $r$-partite graphs where $r\ge 4$. Dai, Liu, and Zhang \cite{DLZ2025} showed that $\delta(G)\ge (r-1 - \frac{1}{2s^2})n$ suffices for a copy of $K_r(s)$ in $G\in \mathcal{G}_r(n)$. 
In this paper we give an asymptotically tight answer to this problem for all even $r$. To be consistent with Theorem~\ref{mainresult}, we state the problem in the complementary form. Given an $r$-partite graph, an \textit{$s$-blowup of independent transversal}, denoted by $IT(s)$, is an independent set containing exactly $s$ vertices from each part. 
	
\begin{thm}\label{thm:blowup}
	For any even integer $r\ge 2$ and integer $s\ge 2$, there exist constants $C_{r,s}$ and $n_0$ such that every graph $G\in \mathcal{G}_r(n)$ with $n\ge n_0$ and maximum degree $\Delta(G)\le \frac{r}{2r-2}n- C_{r,s}n^{1-1/s}$ contains an $IT(s)$. 
\end{thm}

We will show that Theorem~\ref{thm:blowup} is best possible in Section 2, provided that $z(n; s)=\Omega(n^{2- 1/s})$, which is known for $s=2, 3$ \cite{B1966,ERS1966}. 

\medskip
Both proofs of Theorems \ref{mainresult} and \ref{thm:blowup} reply on the Graph Removal Lemma and structure theorems of Haxell and Szab\'o \cite{HS2006} on graphs without an independent transversal. Indeed, to prove Theorem \ref{mainresult}, we first apply the Removal Lemma and the results of \cite{HS2006} to show that every graph $G\in \mathcal{G}_r(n)$ without many independent transversals is not far from the union $G'$ of $r-1$ vertex-disjoint complete bipartite graphs. A key property of $G'$ is that, given any $V_i$ of $G$, every set of $2r-2$ representatives (one vertex from each \emph{partition set} of $G'$ with respect to the $r-1$ complete bipartite graphs that form $G'$) contains a partial independent transversal of size $r-1$ that avoids $V_i$. We then use the maximum degree condition of $G$ to find a vertex not adjacent to these $2r-2$ representatives and thus obtain a full independent transversal of $G$. 
The number of independent transversals can be easily calculated by using the number of the $2r-2$ representatives and the number of the vertices added at last.
The proof of Theorem \ref{thm:blowup} is more complicated because, for example, there may not be $s\ge 2$ vertices from the same $V_i$ that are not adjacent to the $2r-2$ representatives in $G$. Indeed, to prove Theorem \ref{thm:blowup}, we consider the set $V_0=V(G)\setminus V(G')$. If there are many vertices of $V_0$ that have many non-neighbors of $G$ in each partition set of $G'$, then we use double counting arguments to obtain an $IT(s)$. Otherwise, most vertices in $V_0$ are adjacent to most vertices in \emph{some} (unique) partition set of $G'$, and we add these vertices of $V_0$ to $G'$ accordingly.
By the maximum degree of $G$, one of the resulting $r-1$ bipartite graphs must have many non-edges of $G$, thus containing a copy of $K_{s,s}$ in the complement of $G$. We then extend the vertices in this $K_{s,s}$ to an $IT(s)$ of $G$ by using the structure of $G'$. 

\medskip
Let us mention some other results on the existence of independent transversals. Loh and Sudakov \cite{LS2007} proved that every $r$-partite
graph $G$ with parts $V_1, \dots, V_r$ of size $n\ge (1+o(1))D$, $\Delta(G)\le D$, and $\max_{i\ne j} \Delta(G[V_i,V_j])= o(D)$ contains an independent transversal. 
Recently, Glock and Sudakov \cite{GS2020} and independently Kang and Kelly \cite{KK2021} strengthened the result of \cite{LS2007} by replacing $\Delta(G)\le D$ by $b(G):=\max_{i} \sum_{v\in V_i} d(v)/n \le D$. 
All these results were obtained by the semi-random method.
Using a counting argument, Wanless and Wood \cite{WW2022} showed that every $G\in \mathcal{G}_r(n)$ with $b(G)\le n/4$ contains at least $(n/2)^r$ independent transversals (the existence of an independent transversal was noted in \cite{KK2021}). Groenland, Kaiser, Treffers and Wales \cite{GKTW2023} gave a construction with $b(G)= (1+o(1))n/4$ and without any independent transversal.
Haxell and Wdowinski \cite{HW2024JGT} showed that $\Delta(G)\le \alpha n$ and $b(G)\le \beta n$ guarantees an independent transversal if and only if 
$\alpha \le 1/2$ or $\beta \le 1/4$ or $\beta \le 2\alpha(1-\alpha)$.

\subsection{Notation}
Let $G=(V, E)$ be a graph and let $A,B$ be two disjoint subsets of $V$. We write $v(G)=|V|$ and $e(G)=|E|$. We denote by $G[A]$ the induced subgraph of $G$ on $A$, and by $G[A,B]$ the induced bipartite subgraph of $G$ between $A$ and $B$. 
For $x\in V(G)$, we let $N_G(x)$ be the neighborhood of $x$ in $G$, and $N_G(x,A):=N_G(x)\cap A$. 
The degree of $x$ is $d_G(x)= |N_G(x)|$ and the degree of $x$ in $A$ is $d_G(x, A)= |N_G(x, A)|$. We let $\overline{d}_G(x,A):=|A\setminus (N_G(x)\cup\{x\})|$ be the number of nonneighbors of $x$ in $A$.
Given an edge set $E_0$, we let $d_{E_0}(x):=|\{e\in E_0: e\ni x\}|$. 

We also write $N_G(A):=\cup_{x\in A}N_G(x)$. We use $e_G(A,B)$ to denote the number of edges between $A$ and $B$ in $G$, and use $\overline{e}_G(A,B)$ to denote the number of nonadjacent pairs between $A$ and $B$ in $G$ (including pairs of vertices from the same part of $G$ when $G$ is multipartite). 
We omit the subscript in these notation if the underlying graph is clear from the context.

Let $G=(V, E)$ and $G'=(V', E')$ be two graphs. 
We define $G\setminus G':= (V, E\setminus E')$ and $G\cap G':=(V\cap V',E\cap E')$.
We write $G=\emptyset$ when $E=\emptyset$.
For $E_0\subseteq E$, we use $G-E_0$ to denote the graph obtained by removing the edges of $E_0$ from $G$. 
We denote by $\overline{G}$ the complement graph of $G$ whose edge set consists of all nonadjacent pairs of $V$ under $G$. 

Let $K_r$ denote the complete graph on $r$ vertices and $K_r(s)$ denote the complete $r$-partite graph with $s$ vertices in each part. We often write $K_2(s)$ as $K_{s,s}$.
Given disjoint sets $V_1,\cdots, V_r$, we let $K(V_1,V_2,\cdots,V_r)$ be the complete $r$-partite graph with parts $V_1,\cdots,V_r$.

Let $[n]=\{1, \dots, n\}$ for $n\in \mathbb{N}$. We often use $uv$ to abbreviate $\{u, v\}$. 
We write $a\ll b$ if there exists an increasing positive function $f$ such that $a\le f(b)$. We will omit floors and ceilings when these are not crucial.

\subsection{Organization of the paper}
%The organization of the paper is as follows. 
In Section 2 we give a construction that proves the upper bound in Theorem \ref{mainresult} and the tightness of Theorem \ref{thm:blowup}. In Section 3 we introduce necessary tools including the concept of \textit{induced matching configuration} introduced by Haxell and Szab\'o \cite{HS2006}. We prove Theorem \ref{mainresult} (lower bound) in Section 4 and Theorem \ref{thm:blowup} in Section 5, and give concluding remarks in the last section.

%-------------------------------------
\section{Constructions}

In this section, we construct extremal graphs for Theorems \ref{mainresult} and \ref{thm:blowup}, and prove the upper bound in Theorem \ref{mainresult}. 

First let us recall the construction of Szab\'o and Tardos \cite[Construction 3.3]{ST2006} that gives the (tight) upper bound of $\Delta_r(n)$ for any even $r$.

\begin{cons}[\cite{ST2006}]\label{cons0}
    Suppose $r, n, t$ are positive integers such that $r$ is even and $r-1$ divides $n$. Let $A_1,B_1,\cdots, A_{r-1},B_{r-1}$ be disjoint sets of $\frac{rn}{2r-2}$ and let $G_1$ be the union of complete bipartite graphs $K(A_i,B_i)$ for $i=1, \dots, r-1$. Thus, $G_1$ is $\frac{rn}{2r-2}$-regular.
    We partition $V(G_1)=\bigcup_{i=1}^{r-1} A_i\cup B_i$ into parts $V_1, \dots, V_n$ of size $n$ as follows. For $i=1,\cdots,r/2-1$, let $V_i$ be the union of $A_i$ and a subset $B_{i+1}'\subseteq B_{i+1}$ of size $|B_{i+1}'|=n-|A_i|$, and let $V_{r/2}=B_1\cup (\cup_{i=2}^{r/2}B_i'')$ where $B_i''=B_i\setminus B_i'$. For $i=r/2+1,\cdots, r-1$, let $V_i$ be the union of $B_i$ and a subset $A_{i-1}'\subseteq A_{i-1}$ of size $|A_{i-1}'|=n-|B_i|$, and let $V_r=A_{r-1}\cup_{i=r/2}^{r-2}A_i''$ where $A_i''=A_i\setminus A_i'$. 
\end{cons}
It is easy to see that $G_1$ given by Construction~\ref{cons0} contains no independent transversal. Indeed, since $G_1[A_{i}, B_{i}]$ for all $i$ is complete, an independent transversal $T$ must avoid $A_{r/2}$ or $B_{r/2}$. Without loss of generality, assume $T\cap B_{r/2}=\emptyset$. Then the structure of $G$ forces $T\cap A_i\ne\emptyset$ for all $i<r/2$. It follows that $T\cap B_i=\emptyset$ for all $i\le r/2$ and consequently, $T\cap V_{r/2}= \emptyset$ implying that $T$ is not a transversal.

\medskip

To obtain a graph $G\in \mathcal{G}_r(n)$ with $\Delta(G)=\frac{rn}{2r-2}-t$, one can remove a $t$-regular spanning subgraph from $G_1$ (such subgraph exists because $G_1$ is the union of balanced complete bipartite graphs). Since every independent transversal must contain one removed edge, the number of independent transversals is at most $(\frac{trn}{2}) n^{r-2} = \frac{r}{2}tn^{r-1}$ (a more careful case analysis can give a better bound). Below we present a different construction by altering the size of $A_i, B_i$ and removing some edges from $G_1[A_{1}, B_{1}]$ and $G_1[A_{r}, B_{r}]$. This construction allows us to prove the tightness of Theorem \ref{thm:blowup}.

\begin{cons}\label{cons1}
        Suppose $r, n, t$ are positive integers and $r\ge 2$ is even. Let $D=\lceil\frac{rn}{2r-2}\rceil$. Given a balanced bipartite graph $H$ with parts of size $D+\frac{(r-2)t}{2}$,
        we construct a graph $G_2= G_2(t,H)$ with $V(G)= \bigcup_{i=1}^{r-1} A_i\cup B_i = \bigcup_{i=1}^r V_i$ as follows. Let $A_1,B_1,A_{r-1},B_{r-1}$ be sets of size $D+\frac{(r-2)t}{2}$. Let $A_{r/2},B_{r/2}$ be sets of size $\frac{rn}{2}-(r-2)D-2t$ and let $A_i,B_i$, $i\in\{2,\cdots,\frac{r}{2}-1,\frac{r}{2}+1,\cdots, r-2\}$, be sets of size $D-t$. For $i=1,\cdots,r/2-1$, let $V_i$ be the union of $A_i$ and a subset $B_{i+1}'\subseteq B_{i+1}$ of size $|B_{i+1}'|=n-|A_i|$, and let $V_{r/2}=B_1\cup (\cup_{i=2}^{r/2}B_i'')$ where $B_i''=B_i\setminus B_i'$. For $i=r/2+1,\cdots, r-1$, let $V_i$ be the union of $B_i$ and a subset $A_{i-1}'\subseteq A_{i-1}$ of size $|A_{i-1}'|=n-|B_i|$, and let $V_r=A_{r-1}\cup_{i=r/2}^{r-2}A_i''$ where $A_i''=A_i\setminus A_i'$.  
    We let $G_2[A_i,B_i]$ be complete for $2\le i\le r-2$ and let $G_2[A_1,B_1]\cong G_2[A_{r-1},B_{r-1}]\cong K_{D+\frac{(r-2)t}{2}, D+\frac{(r-2)t}{2}}- H$ be the bipartite complement of $H$. (See Figure~\ref{fig_construction}.)
\end{cons} 

When $r=2$, we have $B_1=V_1$ and $A_1=V_2$ with $|A_1|=|B_1|=D=n$, and 
$G_2= K(A_1, B_1)\setminus H$ is the bipartite complement of $H$.

When $r\ge 4$, it is clear that $\Delta(G_2)=\max\{D - t, D + \frac{(r-2)t}{2}- \delta(H)\}$. Let us verify that $|V_i|=n$ for $i\in [r]$. This is clear for $i\ne \frac{r}{2},r$. We know $|V_{r/2}|=|V_r|$ by symmetry. Since 
\begin{align*}
  \sum_{i=1}^r|V_i|&=|V(G)|=\sum_{i=1}^{r-1}\left(|A_i|+|B_i|\right)\\
    &=4\left(D+\frac{(r-2)t}{2}\right)+2(r-4)(D-t)+2\left(\frac{rn}{2}-(r-2)D-2t\right)=rn,  
\end{align*}
it follows that $|V_{r/2}|=|V_r|=n$.

\tikzset{every picture/.style={line width=0.75pt,scale=0.8}} %set default line width to 0.75pt        
\begin{figure}[ht]
\centering
\begin{tikzpicture}[x=0.75pt,y=0.75pt,yscale=-1,xscale=1]

%Rounded Rect [id:dp973834271906207] 
\draw   (113.33,101) .. controls (113.33,96.58) and (116.92,93) .. (121.33,93) -- (175.33,93) .. controls (179.75,93) and (183.33,96.58) .. (183.33,101) -- (183.33,125) .. controls (183.33,129.42) and (179.75,133) .. (175.33,133) -- (121.33,133) .. controls (116.92,133) and (113.33,129.42) .. (113.33,125) -- cycle ;
%Shape: Rectangle [id:dp8464741763595858] 
\draw  [color={rgb, 255:red, 74; green, 144; blue, 226 }  ,draw opacity=0.3 ][fill={rgb, 255:red, 74; green, 144; blue, 226 }  ,fill opacity=0.3 ] (123.33,133) -- (173.33,133) -- (173.33,213) -- (123.33,213) -- cycle ;
%Shape: Boxed Line [id:dp6366135492111213] 
\draw    (148.33,113) ;
\draw [shift={(148.33,113)}, rotate = 0] [color={rgb, 255:red, 0; green, 0; blue, 0 }  ][fill={rgb, 255:red, 0; green, 0; blue, 0 }  ][line width=0.75]      (0, 0) circle [x radius= 3.35, y radius= 3.35]   ;
%Straight Lines [id:da19957378942656945] 
\draw    (173.33,113) ;
\draw [shift={(173.33,113)}, rotate = 0] [color={rgb, 255:red, 0; green, 0; blue, 0 }  ][fill={rgb, 255:red, 0; green, 0; blue, 0 }  ][line width=0.75]      (0, 0) circle [x radius= 3.35, y radius= 3.35]   ;
%Straight Lines [id:da8669699784944238] 
\draw    (123.33,113) ;
\draw [shift={(123.33,113)}, rotate = 0] [color={rgb, 255:red, 0; green, 0; blue, 0 }  ][fill={rgb, 255:red, 0; green, 0; blue, 0 }  ][line width=0.75]      (0, 0) circle [x radius= 3.35, y radius= 3.35]   ;
%Rounded Rect [id:dp8000666062428581] 
\draw   (113.33,221) .. controls (113.33,216.58) and (116.92,213) .. (121.33,213) -- (175.33,213) .. controls (179.75,213) and (183.33,216.58) .. (183.33,221) -- (183.33,245) .. controls (183.33,249.42) and (179.75,253) .. (175.33,253) -- (121.33,253) .. controls (116.92,253) and (113.33,249.42) .. (113.33,245) -- cycle ;
%Shape: Boxed Line [id:dp06819391266199992] 
\draw    (148.33,233) ;
\draw [shift={(148.33,233)}, rotate = 0] [color={rgb, 255:red, 0; green, 0; blue, 0 }  ][fill={rgb, 255:red, 0; green, 0; blue, 0 }  ][line width=0.75]      (0, 0) circle [x radius= 3.35, y radius= 3.35]   ;
%Straight Lines [id:da9832304893362991] 
\draw    (173.33,233) ;
\draw [shift={(173.33,233)}, rotate = 0] [color={rgb, 255:red, 0; green, 0; blue, 0 }  ][fill={rgb, 255:red, 0; green, 0; blue, 0 }  ][line width=0.75]      (0, 0) circle [x radius= 3.35, y radius= 3.35]   ;
%Straight Lines [id:da5039080883175142] 
\draw    (123.33,233) ;
\draw [shift={(123.33,233)}, rotate = 0] [color={rgb, 255:red, 0; green, 0; blue, 0 }  ][fill={rgb, 255:red, 0; green, 0; blue, 0 }  ][line width=0.75]      (0, 0) circle [x radius= 3.35, y radius= 3.35]   ;
%Rounded Rect [id:dp8564801970325435] 
\draw   (223.33,101) .. controls (223.33,96.58) and (226.92,93) .. (231.33,93) -- (285.33,93) .. controls (289.75,93) and (293.33,96.58) .. (293.33,101) -- (293.33,125) .. controls (293.33,129.42) and (289.75,133) .. (285.33,133) -- (231.33,133) .. controls (226.92,133) and (223.33,129.42) .. (223.33,125) -- cycle ;
%Shape: Rectangle [id:dp2984539373519557] 
\draw  [color={rgb, 255:red, 74; green, 144; blue, 226 }  ,draw opacity=0.3 ][fill={rgb, 255:red, 74; green, 144; blue, 226 }  ,fill opacity=0.3 ] (233.33,133) -- (283.33,133) -- (283.33,213) -- (233.33,213) -- cycle ;
%Shape: Boxed Line [id:dp5171618551280572] 
\draw    (258.33,113) ;
\draw [shift={(258.33,113)}, rotate = 0] [color={rgb, 255:red, 0; green, 0; blue, 0 }  ][fill={rgb, 255:red, 0; green, 0; blue, 0 }  ][line width=0.75]      (0, 0) circle [x radius= 3.35, y radius= 3.35]   ;
%Straight Lines [id:da17470525425390715] 
\draw    (283.33,113) ;
\draw [shift={(283.33,113)}, rotate = 0] [color={rgb, 255:red, 0; green, 0; blue, 0 }  ][fill={rgb, 255:red, 0; green, 0; blue, 0 }  ][line width=0.75]      (0, 0) circle [x radius= 3.35, y radius= 3.35]   ;
%Straight Lines [id:da026558361047256174] 
\draw    (233.33,113) ;
\draw [shift={(233.33,113)}, rotate = 0] [color={rgb, 255:red, 0; green, 0; blue, 0 }  ][fill={rgb, 255:red, 0; green, 0; blue, 0 }  ][line width=0.75]      (0, 0) circle [x radius= 3.35, y radius= 3.35]   ;
%Rounded Rect [id:dp8127689335459987] 
\draw   (223.33,221) .. controls (223.33,216.58) and (226.92,213) .. (231.33,213) -- (285.33,213) .. controls (289.75,213) and (293.33,216.58) .. (293.33,221) -- (293.33,245) .. controls (293.33,249.42) and (289.75,253) .. (285.33,253) -- (231.33,253) .. controls (226.92,253) and (223.33,249.42) .. (223.33,245) -- cycle ;
%Shape: Boxed Line [id:dp8599335006327484] 
\draw    (258.33,233) ;
\draw [shift={(258.33,233)}, rotate = 0] [color={rgb, 255:red, 0; green, 0; blue, 0 }  ][fill={rgb, 255:red, 0; green, 0; blue, 0 }  ][line width=0.75]      (0, 0) circle [x radius= 3.35, y radius= 3.35]   ;
%Straight Lines [id:da01541496057916314] 
\draw    (283.33,233) ;
\draw [shift={(283.33,233)}, rotate = 0] [color={rgb, 255:red, 0; green, 0; blue, 0 }  ][fill={rgb, 255:red, 0; green, 0; blue, 0 }  ][line width=0.75]      (0, 0) circle [x radius= 3.35, y radius= 3.35]   ;
%Straight Lines [id:da7818927876595139] 
\draw    (233.33,233) ;
\draw [shift={(233.33,233)}, rotate = 0] [color={rgb, 255:red, 0; green, 0; blue, 0 }  ][fill={rgb, 255:red, 0; green, 0; blue, 0 }  ][line width=0.75]      (0, 0) circle [x radius= 3.35, y radius= 3.35]   ;
%Rounded Rect [id:dp9713983775398669] 
\draw   (333.33,101) .. controls (333.33,96.58) and (336.92,93) .. (341.33,93) -- (395.33,93) .. controls (399.75,93) and (403.33,96.58) .. (403.33,101) -- (403.33,125) .. controls (403.33,129.42) and (399.75,133) .. (395.33,133) -- (341.33,133) .. controls (336.92,133) and (333.33,129.42) .. (333.33,125) -- cycle ;
%Shape: Rectangle [id:dp5719503922973632] 
\draw  [color={rgb, 255:red, 74; green, 144; blue, 226 }  ,draw opacity=0.3 ][fill={rgb, 255:red, 74; green, 144; blue, 226 }  ,fill opacity=0.3 ] (343.33,133) -- (393.33,133) -- (393.33,213) -- (343.33,213) -- cycle ;
%Shape: Boxed Line [id:dp7891478095032562] 
\draw    (368.33,113) ;
\draw [shift={(368.33,113)}, rotate = 0] [color={rgb, 255:red, 0; green, 0; blue, 0 }  ][fill={rgb, 255:red, 0; green, 0; blue, 0 }  ][line width=0.75]      (0, 0) circle [x radius= 3.35, y radius= 3.35]   ;
%Straight Lines [id:da20155744811865484] 
\draw    (393.33,113) ;
\draw [shift={(393.33,113)}, rotate = 0] [color={rgb, 255:red, 0; green, 0; blue, 0 }  ][fill={rgb, 255:red, 0; green, 0; blue, 0 }  ][line width=0.75]      (0, 0) circle [x radius= 3.35, y radius= 3.35]   ;
%Straight Lines [id:da48453325293474125] 
\draw    (343.33,113) ;
\draw [shift={(343.33,113)}, rotate = 0] [color={rgb, 255:red, 0; green, 0; blue, 0 }  ][fill={rgb, 255:red, 0; green, 0; blue, 0 }  ][line width=0.75]      (0, 0) circle [x radius= 3.35, y radius= 3.35]   ;
%Rounded Rect [id:dp7073353503379989] 
\draw   (333.33,221) .. controls (333.33,216.58) and (336.92,213) .. (341.33,213) -- (395.33,213) .. controls (399.75,213) and (403.33,216.58) .. (403.33,221) -- (403.33,245) .. controls (403.33,249.42) and (399.75,253) .. (395.33,253) -- (341.33,253) .. controls (336.92,253) and (333.33,249.42) .. (333.33,245) -- cycle ;
%Shape: Boxed Line [id:dp10033974829647696] 
\draw    (368.33,233) ;
\draw [shift={(368.33,233)}, rotate = 0] [color={rgb, 255:red, 0; green, 0; blue, 0 }  ][fill={rgb, 255:red, 0; green, 0; blue, 0 }  ][line width=0.75]      (0, 0) circle [x radius= 3.35, y radius= 3.35]   ;
%Straight Lines [id:da9638750974497918] 
\draw    (393.33,233) ;
\draw [shift={(393.33,233)}, rotate = 0] [color={rgb, 255:red, 0; green, 0; blue, 0 }  ][fill={rgb, 255:red, 0; green, 0; blue, 0 }  ][line width=0.75]      (0, 0) circle [x radius= 3.35, y radius= 3.35]   ;
%Straight Lines [id:da7451162719241766] 
\draw    (343.33,233) ;
\draw [shift={(343.33,233)}, rotate = 0] [color={rgb, 255:red, 0; green, 0; blue, 0 }  ][fill={rgb, 255:red, 0; green, 0; blue, 0 }  ][line width=0.75]      (0, 0) circle [x radius= 3.35, y radius= 3.35]   ;
%Rounded Rect [id:dp9080440551238329] 
\draw   (443.33,101) .. controls (443.33,96.58) and (446.92,93) .. (451.33,93) -- (505.33,93) .. controls (509.75,93) and (513.33,96.58) .. (513.33,101) -- (513.33,125) .. controls (513.33,129.42) and (509.75,133) .. (505.33,133) -- (451.33,133) .. controls (446.92,133) and (443.33,129.42) .. (443.33,125) -- cycle ;
%Shape: Rectangle [id:dp0644688933593911] 
\draw  [color={rgb, 255:red, 74; green, 144; blue, 226 }  ,draw opacity=0.3 ][fill={rgb, 255:red, 74; green, 144; blue, 226 }  ,fill opacity=0.3 ] (453.33,133) -- (503.33,133) -- (503.33,213) -- (453.33,213) -- cycle ;
%Shape: Boxed Line [id:dp28489168738143755] 
\draw    (478.33,113) ;
\draw [shift={(478.33,113)}, rotate = 0] [color={rgb, 255:red, 0; green, 0; blue, 0 }  ][fill={rgb, 255:red, 0; green, 0; blue, 0 }  ][line width=0.75]      (0, 0) circle [x radius= 3.35, y radius= 3.35]   ;
%Straight Lines [id:da43882586444926797] 
\draw    (503.33,113) ;
\draw [shift={(503.33,113)}, rotate = 0] [color={rgb, 255:red, 0; green, 0; blue, 0 }  ][fill={rgb, 255:red, 0; green, 0; blue, 0 }  ][line width=0.75]      (0, 0) circle [x radius= 3.35, y radius= 3.35]   ;
%Straight Lines [id:da4338373354251519] 
\draw    (453.33,113) ;
\draw [shift={(453.33,113)}, rotate = 0] [color={rgb, 255:red, 0; green, 0; blue, 0 }  ][fill={rgb, 255:red, 0; green, 0; blue, 0 }  ][line width=0.75]      (0, 0) circle [x radius= 3.35, y radius= 3.35]   ;
%Rounded Rect [id:dp9707641489128669] 
\draw   (443.33,221) .. controls (443.33,216.58) and (446.92,213) .. (451.33,213) -- (505.33,213) .. controls (509.75,213) and (513.33,216.58) .. (513.33,221) -- (513.33,245) .. controls (513.33,249.42) and (509.75,253) .. (505.33,253) -- (451.33,253) .. controls (446.92,253) and (443.33,249.42) .. (443.33,245) -- cycle ;
%Shape: Boxed Line [id:dp562578085594144] 
\draw    (478.33,233) ;
\draw [shift={(478.33,233)}, rotate = 0] [color={rgb, 255:red, 0; green, 0; blue, 0 }  ][fill={rgb, 255:red, 0; green, 0; blue, 0 }  ][line width=0.75]      (0, 0) circle [x radius= 3.35, y radius= 3.35]   ;
%Straight Lines [id:da6644726705784771] 
\draw    (503.33,233) ;
\draw [shift={(503.33,233)}, rotate = 0] [color={rgb, 255:red, 0; green, 0; blue, 0 }  ][fill={rgb, 255:red, 0; green, 0; blue, 0 }  ][line width=0.75]      (0, 0) circle [x radius= 3.35, y radius= 3.35]   ;
%Straight Lines [id:da7689779991130714] 
\draw    (453.33,233) ;
\draw [shift={(453.33,233)}, rotate = 0] [color={rgb, 255:red, 0; green, 0; blue, 0 }  ][fill={rgb, 255:red, 0; green, 0; blue, 0 }  ][line width=0.75]      (0, 0) circle [x radius= 3.35, y radius= 3.35]   ;
%Rounded Rect [id:dp7339666882259515] 
\draw   (553.33,101) .. controls (553.33,96.58) and (556.92,93) .. (561.33,93) -- (615.33,93) .. controls (619.75,93) and (623.33,96.58) .. (623.33,101) -- (623.33,125) .. controls (623.33,129.42) and (619.75,133) .. (615.33,133) -- (561.33,133) .. controls (556.92,133) and (553.33,129.42) .. (553.33,125) -- cycle ;
%Shape: Rectangle [id:dp4479355875325266] 
\draw  [color={rgb, 255:red, 74; green, 144; blue, 226 }  ,draw opacity=0.3 ][fill={rgb, 255:red, 74; green, 144; blue, 226 }  ,fill opacity=0.3 ] (563.33,133) -- (613.33,133) -- (613.33,213) -- (563.33,213) -- cycle ;
%Shape: Boxed Line [id:dp1274423496483823] 
\draw    (588.33,113) ;
\draw [shift={(588.33,113)}, rotate = 0] [color={rgb, 255:red, 0; green, 0; blue, 0 }  ][fill={rgb, 255:red, 0; green, 0; blue, 0 }  ][line width=0.75]      (0, 0) circle [x radius= 3.35, y radius= 3.35]   ;
%Straight Lines [id:da9716276725363966] 
\draw    (613.33,113) ;
\draw [shift={(613.33,113)}, rotate = 0] [color={rgb, 255:red, 0; green, 0; blue, 0 }  ][fill={rgb, 255:red, 0; green, 0; blue, 0 }  ][line width=0.75]      (0, 0) circle [x radius= 3.35, y radius= 3.35]   ;
%Straight Lines [id:da5561969683005266] 
\draw    (563.33,113) ;
\draw [shift={(563.33,113)}, rotate = 0] [color={rgb, 255:red, 0; green, 0; blue, 0 }  ][fill={rgb, 255:red, 0; green, 0; blue, 0 }  ][line width=0.75]      (0, 0) circle [x radius= 3.35, y radius= 3.35]   ;
%Rounded Rect [id:dp5270796940620905] 
\draw   (553.33,221) .. controls (553.33,216.58) and (556.92,213) .. (561.33,213) -- (615.33,213) .. controls (619.75,213) and (623.33,216.58) .. (623.33,221) -- (623.33,245) .. controls (623.33,249.42) and (619.75,253) .. (615.33,253) -- (561.33,253) .. controls (556.92,253) and (553.33,249.42) .. (553.33,245) -- cycle ;
%Shape: Boxed Line [id:dp27336005003364283] 
\draw    (588.33,233) ;
\draw [shift={(588.33,233)}, rotate = 0] [color={rgb, 255:red, 0; green, 0; blue, 0 }  ][fill={rgb, 255:red, 0; green, 0; blue, 0 }  ][line width=0.75]      (0, 0) circle [x radius= 3.35, y radius= 3.35]   ;
%Straight Lines [id:da20926844632527541] 
\draw    (613.33,233) ;
\draw [shift={(613.33,233)}, rotate = 0] [color={rgb, 255:red, 0; green, 0; blue, 0 }  ][fill={rgb, 255:red, 0; green, 0; blue, 0 }  ][line width=0.75]      (0, 0) circle [x radius= 3.35, y radius= 3.35]   ;
%Straight Lines [id:da5134104429863886] 
\draw    (563.33,233) ;
\draw [shift={(563.33,233)}, rotate = 0] [color={rgb, 255:red, 0; green, 0; blue, 0 }  ][fill={rgb, 255:red, 0; green, 0; blue, 0 }  ][line width=0.75]      (0, 0) circle [x radius= 3.35, y radius= 3.35]   ;
%Shape: Polygon Curved [id:ds6664905680285884] 
\draw  [fill={rgb, 255:red, 155; green, 155; blue, 155 }  ,fill opacity=0.15 ][dash pattern={on 4.5pt off 4.5pt}] (113.11,88.44) .. controls (118.44,79.78) and (169.11,77.78) .. (179.78,87.11) .. controls (190.44,96.44) and (273.11,221.11) .. (267.78,235.11) .. controls (262.44,249.11) and (243.11,247.78) .. (231.78,242.44) .. controls (220.44,237.11) and (111.11,140.44) .. (108.44,128.44) .. controls (105.78,116.44) and (107.78,97.11) .. (113.11,88.44) -- cycle ;
%Shape: Polygon Curved [id:ds8668996040301584] 
\draw  [fill={rgb, 255:red, 155; green, 155; blue, 155 }  ,fill opacity=0.15 ][dash pattern={on 4.5pt off 4.5pt}] (223.44,88.11) .. controls (228.78,79.44) and (279.44,77.44) .. (290.11,86.78) .. controls (300.78,96.11) and (383.44,220.78) .. (378.11,234.78) .. controls (372.78,248.78) and (353.44,247.44) .. (342.11,242.11) .. controls (330.78,236.78) and (221.44,140.11) .. (218.78,128.11) .. controls (216.11,116.11) and (218.11,96.78) .. (223.44,88.11) -- cycle ;
%Shape: Polygon Curved [id:ds31771435130122017] 
\draw  [fill={rgb, 255:red, 155; green, 155; blue, 155 }  ,fill opacity=0.15 ][dash pattern={on 4.5pt off 4.5pt}] (472.11,106.78) .. controls (480.11,90.11) and (506.78,98.11) .. (516.11,104.11) .. controls (525.44,110.11) and (610.11,198.11) .. (624.11,212.11) .. controls (638.11,226.11) and (636.11,248.78) .. (633.44,258.11) .. controls (630.78,267.44) and (570.11,278.78) .. (550.11,248.78) .. controls (530.11,218.78) and (464.11,123.44) .. (472.11,106.78) -- cycle ;
%Shape: Polygon Curved [id:ds17480619273695663] 
\draw  [fill={rgb, 255:red, 155; green, 155; blue, 155 }  ,fill opacity=0.15 ][dash pattern={on 4.5pt off 4.5pt}] (361.44,106.78) .. controls (369.44,90.11) and (396.11,98.11) .. (405.44,104.11) .. controls (414.78,110.11) and (499.44,198.11) .. (513.44,212.11) .. controls (527.44,226.11) and (525.44,248.78) .. (522.78,258.11) .. controls (520.11,267.44) and (459.44,278.78) .. (439.44,248.78) .. controls (419.44,218.78) and (353.44,123.44) .. (361.44,106.78) -- cycle ;
%Shape: Polygon Curved [id:ds9687905933664407] 
\draw  [fill={rgb, 255:red, 155; green, 155; blue, 155 }  ,fill opacity=0.15 ][dash pattern={on 4.5pt off 4.5pt}] (127.44,280.78) .. controls (108.11,269.44) and (94.78,252.11) .. (112.78,225.44) .. controls (130.78,198.78) and (166.56,203.89) .. (183.33,221) .. controls (200.11,238.11) and (209.44,270.11) .. (242.78,272.11) .. controls (276.11,274.11) and (260.11,229.44) .. (284.11,224.78) .. controls (308.11,220.11) and (327.44,278.39) .. (350.11,272.11) .. controls (372.78,265.83) and (378.78,224.78) .. (396.78,222.11) .. controls (414.78,219.44) and (418.11,272.78) .. (393.44,283.44) .. controls (368.78,294.11) and (146.78,292.11) .. (127.44,280.78) -- cycle ;
%Shape: Polygon Curved [id:ds9259106189904314] 
\draw  [fill={rgb, 255:red, 155; green, 155; blue, 155 }  ,fill opacity=0.15 ][dash pattern={on 4.5pt off 4.5pt}] (349.44,67.44) .. controls (366.11,54.11) and (590.78,50.78) .. (613.44,68.11) .. controls (636.11,85.44) and (636.78,108.11) .. (626.11,121.44) .. controls (615.44,134.78) and (570.11,132.11) .. (550.78,114.78) .. controls (531.44,97.44) and (507.44,73.44) .. (488.78,80.78) .. controls (470.11,88.11) and (466.11,123.44) .. (450.78,124.11) .. controls (435.44,124.78) and (407.44,85.44) .. (383.44,83.44) .. controls (359.44,81.44) and (356.11,127.44) .. (340.11,123.44) .. controls (324.11,119.44) and (332.78,80.78) .. (349.44,67.44) -- cycle ;
%Straight Lines [id:da31441481158579254] 
\draw [color={rgb, 255:red, 208; green, 2; blue, 27 }  ,draw opacity=1 ] [dash pattern={on 4.5pt off 4.5pt}]  (136,121.33) -- (135.33,227.33) ;
%Straight Lines [id:da3598789800568487] 
\draw [color={rgb, 255:red, 208; green, 2; blue, 27 }  ,draw opacity=1 ] [dash pattern={on 4.5pt off 4.5pt}]  (162.67,120) -- (162,226) ;
%Straight Lines [id:da9706639616902175] 
\draw [color={rgb, 255:red, 208; green, 2; blue, 27 }  ,draw opacity=1 ] [dash pattern={on 4.5pt off 4.5pt}]  (578,121.33) -- (577.33,227.33) ;
%Straight Lines [id:da03527038462021337] 
\draw [color={rgb, 255:red, 208; green, 2; blue, 27 }  ,draw opacity=1 ] [dash pattern={on 4.5pt off 4.5pt}]  (603.33,120) -- (602.67,226) ;

% Text Node
\draw (185.67,158.07) node [anchor=north west][inner sep=0.75pt]    {$V_{1}$};
% Text Node
\draw (297,159.07) node [anchor=north west][inner sep=0.75pt]    {$V_{2}$};
% Text Node
\draw (184.33,263.73) node [anchor=north west][inner sep=0.75pt]    {$V_{3}$};
% Text Node
\draw (421,160.4) node [anchor=north west][inner sep=0.75pt]    {$V_{4}$};
% Text Node
\draw (523,161.4) node [anchor=north west][inner sep=0.75pt]    {$V_{5}$};
% Text Node
\draw (537.67,70.07) node [anchor=north west][inner sep=0.75pt]    {$V_{6}$};
% Text Node
\draw (81,101.4) node [anchor=north west][inner sep=0.75pt]    {$A_{1}$};
% Text Node
\draw (635.67,103.07) node [anchor=north west][inner sep=0.75pt]    {$A_{5}$};
% Text Node
\draw (83,225.4) node [anchor=north west][inner sep=0.75pt]    {$B_{1}$};
% Text Node
\draw (641.67,224.4) node [anchor=north west][inner sep=0.75pt]    {$B_{5}$};
% Text Node
\draw (103.33,38.73) node [anchor=north west][inner sep=0.75pt]    {$\lceil\frac{3n}{5}\rceil+ 2t $};
% Text Node
\draw (229,40.07) node [anchor=north west][inner sep=0.75pt]    {$\lceil\frac{3n}{5}\rceil-t$};
% Text Node
\draw (547,38.4) node [anchor=north west][inner sep=0.75pt]    {$\lceil\frac{3n}{5}\rceil+2t$};
% Text Node
\draw (310,41.07) node [anchor=north west][inner sep=0.75pt]    {$3n-4\lceil\frac{3n}{5}\rceil-2t$};
% Text Node
\draw (447,40.07) node [anchor=north west][inner sep=0.75pt]    {$\lceil\frac{3n}{5}\rceil-t$};
% Text Node
\draw (80,160) node [anchor=north west][inner sep=0.75pt]    {$\textcolor{red}{H}$};
% Text Node
\draw (630,160) node [anchor=north west][inner sep=0.75pt]    {$\textcolor{red}{H}$};

\end{tikzpicture}
\caption{Case $r=6$.}
\label{fig_construction}
\end{figure}
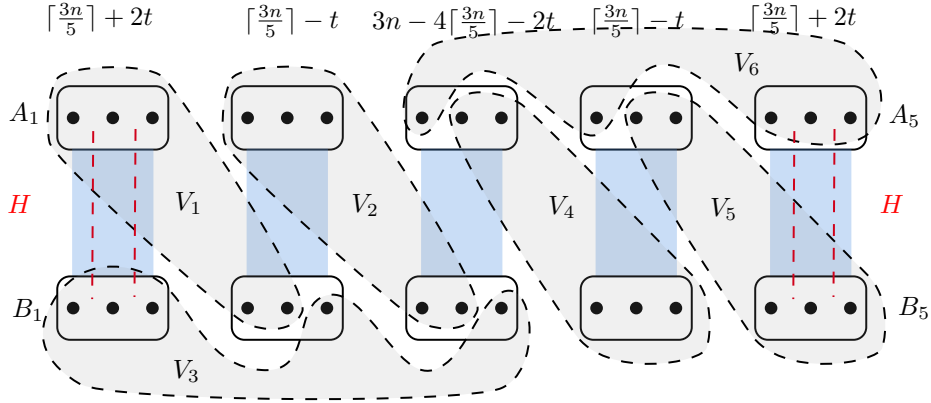

\medskip
The following proposition proves the upper bound in Theorem \ref{mainresult}. 

\begin{prop}\label{pro52}
    Suppose $r, n, t$ are positive integers and $r\ge 4$ is even.
	Let $G_2$ be the graph as in Construction~\ref{cons1} while $H$ is a $\frac{rt}{2} $-regular bipartite graph with parts of size $\lceil\frac{rn}{2r-2}\rceil+\frac{(r-2)t}{2}$. Then $\Delta(G)=\lceil\frac{rn}{2r-2}\rceil-t$, and the number of independent transversals in $G$ is at most $r^2 t(\frac{rn}{2r-2})^{r-1}$ if $t\le \frac{n}{r-1}\left(\sqrt{\frac{r}{r-1}}-1\right)$.
\end{prop}

\begin{proof}
    Let $D=\lceil\frac{rn}{2r-2}\rceil$. By Construction~\ref{cons1}, 
    the vertices of $A_1\cup B_1\cup A_{r-1}\cup B_{r-1}$ has degree $D+\frac{(r-2)t}{2} - \frac{rt}{2} = D-t$ 
    while the other vertices of $G$ has degree $D-t$ or $\frac{rn}{2}-(r-2)D-2t\le D-2t$. Hence $\Delta(G)\le D - t=\lceil\frac{rn}{2r-2}\rceil-t$.
    
	Let $T=\{v_1,\cdots, v_r\}$ be an independent transversal of $G$, where $v_i\in V_i,i\in [r]$.   
	Since $G[A_{r/2},B_{r/2}]$ is complete, at least one of its two parts omits $T$. Below we assume $T\cap B_{r/2}=\emptyset$ and count the number of choices of $v_1,\cdots, v_r$ in this case.
	
    Since $V_{r/2-1}= A_{r/2-1}\cup B_{r/2}'$, it follows that $v_{r/2-1}\in A_{r/2-1}$.
    Since $G[A_{r/2-1},B_{r/2-1}]$ is complete, it follows that $T\cap B_{r/2-1}=\emptyset$, 
    and consequently, $v_{r/2-2}\in A_{r/2-2}$ (recall that $V_{r/2-2}= A_{r/2-2}\cup B_{r/2-1}'$). Repeating this argument, we derive that $v_{i}\in A_{i}$ for $i=1,\cdots, r/2-1$. This implies that $T\cap (\cup_{i=2}^{r/2}B_i)=\emptyset$, and consequently, $v_{r/2}\in B_1$ (recall that $V_{r/2}=B_1\cup (\cup_{i=2}^{r/2}B_i'')$). 
    Since $v_1\in A_1,v_{r/2}\in B_1$ are non-adjacent, once $v_1$ is fixed, there are at most $\frac{rt}{2}$ choices for $v_{r/2}$ because $H$ is $\frac{rt}{2} $-regular. Together there are at most $\prod_{i=1}^{r/2-1} |A_i|\cdot \frac{rt}{2} $ choices of $v_1,\cdots, v_{r/2}$.
		
	Let $q$ be the smallest $i\ge r/2$ such that $A'_i\cap T = \emptyset$ ($q=r-1$ if no such $i$ exists). By definition, we have $A'_i\cap T\ne \emptyset$ for $r/2 \le i\le q-1$. Since $A'_q\cap T = \emptyset$, the arguments in the previous paragraph show that $T\cap B_j\ne \emptyset$ for $j\ge q+1$.  
    Together we have at most 
    $(\prod_{i=r/2}^{q-1}|A_i'|)\cdot(\prod_{j=q+1}^{r-1}|B_i|)\cdot |V_r|$
    choices of $v_{r/2+1},\cdots, v_r$ (we do not have control on the location of $v_r$).

    Putting these together, the number $M$ of choices of $v_1,\cdots,v_r$ satisfies	
	\begin{align*}
	M \le &\prod_{i=1}^{r/2-1} |A_i|\frac{rt}{2}\cdot \sum_{q=r/2}^{r-1} \left(\prod_{i=r/2}^{q-1}|A_i'|\right)\left(\prod_{i=q+1}^{r-1}|B_i|\right)|V_r|  \\
		\le &\frac{rt}{2}\left(D+\frac{(r-2)t}{2}\right) \left(D-t \right)^{r/2-2} \frac{r}{2} \left(D-t \right)^{r/2-2}\left(D+\frac{(r-2)t}{2}\right) n \\
        = & \frac{r^2}{4} tn \left(D+\frac{(r-2)t}{2}\right)^2 \left(D-t \right)^{r-4}. 	
	\end{align*}
    Note that $D-t=\lceil\frac{rn}{2r-2}\rceil-t\le \frac{rn}{2r-2}+1-t\le \frac{rn}{2r-2}$, and thus
	\[
    \left(D+\frac{(r-2)t} {2}\right)^2
    \le \left(\frac{rn}{2r-2}+\frac{rt}{2}\right)^2 \le \left(\frac{rn}{2r-2}\right)^2 \frac{r}{r-1}
    \]
    as $t\le \frac{n}{r-1}\left(\sqrt{\frac{r}{r-1}}-1\right)$. Consequently, 
    \[
    M \le \frac{r^2}{4} tn \left(\frac{rn}{2r-2}\right)^{r-2} \left(\frac{r}{r-1}\right) = \frac{r^2}{2} t\left(\frac{rn}{2r-2}\right)^{r-1}.
    \]
    
    Since the same bound holds when $T\cap A_{r/2}=\emptyset$, the total number of independent transversals in $G$ is at most $r^2t(\frac{rn}{2r-2})^{r-1}$.
\end{proof}

The following proposition shows that Theorem \ref{thm:blowup} is best possible up to the value of $C_{r,s}$, provided that $z(m;s)= \Omega(m^{2-1/s})$. Since a $K_{s, s}$-free graph $G\in \mathcal{G}_2(m)$ with $e(G)= \Omega(m^{2 - 1/s})$ may have vertices of small degree, we first show how to convert $G$ to a $K_{s, s}$-free graph with large minimum degree.

\begin{prop}\label{lem:lb2}
     Suppose $r\ge 2, s\ge 2$ are integers and $r$ is even.
     If there are positive constants $c, C$ such that $cm^{2-1/s}\le z(m;s)\le Cm^{2-1/s}$ for all sufficiently large $m$, then the following holds.
     
     (1) For sufficiently large $N$, there exists a $K_{s,s}$-free bipartite graph $H$ with parts of size $N$ and minimum degree $\delta(H)\ge \frac{c}{8}(\frac{c}{2C})^{\frac{1}{s-1}}N^{1-1/s}$. 
     
     (2) Suppose $\alpha=\frac{c}{8r}(\frac{c}{2C})^{\frac{1}{s-1}}$,  $t=\lfloor\alpha n^{1-1/s}\rfloor$, and $n$ is sufficiently large. 
     Let  $H$ be a $K_{s,s}$-free bipartite graph with parts of size $N=\lceil\frac{rn}{2r-2}\rceil+\frac{r-2}{2}t$ given by Part~(1). Then the graph $G_2=G_2(t, H)$ given by Construction~\ref{cons1} has maximum degree $\Delta(G)\le \lceil\frac{rn}{2r-2}\rceil-t$ and contains no $s$-blowup of independent transversal.
\end{prop}

\begin{proof}
	(1) Let $\gamma=(\frac{c}{2C})^{\frac{s}{s-1}}$ and $m=N/\gamma$. Let $G$ be a $K_{s,s}$-free bipartite graph with parts of size $m$ and with at least $cm^{2-1/s}$ edges. By a well-known fact, $G$ contains a subgraph $G'$ with minimum degree $\delta(G')\ge d(G)/2\ge \frac{c}{2}m^{1-1/s}$, where $d(G)$ is the average degree of $G$. Trivially $G'$ is a $K_{s,s}$-free bipartite graph. Let $A,B$ be two parts of $G'$. Take a copy $G''$ of $G'$ and suppose that the corresponding parts are $A',B'$. Let $\tilde{G}$ be a disjoint union of $G'$ and $G''$ with parts $A\cup B'$ and $A'\cup B$. Then $\tilde{G}$ is a balanced $K_{s,s}$-free bipartite graph. Let $l$ be size of parts of $\tilde{G}$. Then $e(\tilde{G})\le z(l;s)\le Cl^{2-1/s}$ because $\tilde{G}$ is $K_{s,s}$-free. On the other hand, $e(\tilde{G})\ge \delta(\tilde{G})\cdot l$. Thus we have $\delta(\tilde{G})\le Cl^{1-1/s}$. Since $\delta(\tilde{G})=\delta(G')\ge \frac{c}{2}m^{1-1/s}$, it follows that $l\ge \gamma m$, and $l \ge N$. Additionally, we have $l\le 2m= 2N/\gamma$.
	
	Let $p=N/l$, and thus, $\gamma/2\le p \le 1$. It is well-known that by the probabilistic method, we can obtain a balanced bipartite subgraph $H$ of $\tilde{G}$ with parts of size $N=pl$, and minimum degree 
    \[
    \delta(H)\ge \delta(\tilde{G})\, \frac{p}{2} \ge \frac{c}{2}m^{1-1/s}\cdot \frac{\gamma/2}{2}= \frac{c}{8}\left(\frac{N}{\gamma}\right)^{1-1/s}\gamma=\frac{c}{8}\left(\frac{c}{2C}\right)^{\frac{1}{s-1}}N^{1-1/s}.
    \]

    \noindent(2) By the definitions of $\alpha$ and $t$, the graph $H$ has minimum degree 
    \[
    \delta(H)\ge \frac{c}{8}\left(\frac{c}{2C}\right)^{\frac{1}{s-1}}N^{1-1/s} \ge \frac{c}{8}\left(\frac{c}{2C}\right)^{\frac{1}{s-1}}\left(\frac{n}{2}\right)^{1-1/s} =\frac{r\alpha n^{1-1/s}}{2^{1-1/s}}\ge \frac{rt}{2}.
    \]
    Let $D=\lceil\frac{rn}{2r-2}\rceil$. 
    By Construction~\ref{cons1}, the vertices of $A_1\cup B_1\cup A_{r-1}\cup B_{r-1}$ have degree at most $D+\frac{(r-2)t}{2} - \frac{rt}{2} = D-t$ 
    while the other vertices of $G$ has degree at most $D-t$. Thus, $\Delta(G) \le D-t$.

	Suppose $G$ contains an $s$-blowup $T$ of independent transversal. Then $|T\cap V_i|=s$ for $i\in [r]$. Since $G[A_{r/2},B_{r/2}]$ is complete, we have $T\cap A_{r/2}= \emptyset$ or  $T\cap B_{r/2}= \emptyset$. Without loss of generality, assume that $T\cap B_{r/2}=\emptyset$.
    Then $T\cap V_{r/2-1}\subseteq A_{r/2-1}$ as $V_{r/2-1}= A_{r/2-1}\cup B_{r/2}'$. Since $G[A_{r/2-1},B_{r/2-1}]$ is complete, it follows that $T\cap B_{r/2-1}=\emptyset$, and consequently, $T\cap V_{r/2-2}\subseteq A_{r/2-2}$ as $V_{r/2-2}= A_{r/2-2}\cup B_{r/2-1}'$. Continuing these arguments, we conclude that $T\cap B_i=\emptyset$ and $T\cap V_{i-1}\subseteq A_{i-1}$ for $i=2,\cdots, r/2$. Consequently, $T\cap V_{r/2}\subseteq B_1$ as $V_{r/2}=B_1\cup (\cup_{i=2}^{r/2}B_i'')$. Then $T_1:=T\cap V_1$ and $T_2:=T\cap V_{r/2}$ are two $s$-element sets 
    in $A_1$ and $B_1$, respectively, with no edge in $G[T_1,T_2]$. Thus, $H[T_1,T_2]$ is complete, which contradicts $H$ being $K_{s,s}$-free.
    
\end{proof}

%------------------------------------------
\section{Preliminaries}

One of our main tools is the Graph Removal Lemma \cite{EFR1986,F1995,RS1978}.

\begin{lm}[Graph Removal Lemma]\label{removal}
	For an integer $l\ge 2$ and constant $\mu >0$ there exist $\zeta=\zeta(l,\mu)>0$ and $n_0=n_0(l,\mu)$ so that the following holds. Suppose $F$ is a graph on $l$ vertices and $G$ is a graph on $n\ge n_0$ vertices. If $G$ contains at most $\zeta n^l$ copies of $F$, then one can delete $\mu n^2$ edges of $G$ to make it $F$-free.
\end{lm}

We will use several results in \cite{HS2006} that describe the structures of $r$-partite graphs without independent transversals. Note that these results apply to all $r\ge 2$ (even or odd).

Let $G$ be a graph with parts $V_1,\cdots, V_r$ with $r\ge 2$. A set of vertices $I$ is called an \textit{induced matching configuration} (IMC), if $G[I]$ is a perfect matching and it becomes a tree if contracting the vertices of $I\cap V_i$ into one vertex $v_i$ for all $i\in [r]$. For any $q\in [r]$, we say $T\subseteq V(G)$ is a \textit{$V_q$-avoiding partial independent transversal} (PIT) if $T$ is an independent set of size $r-1$ and $|T\cap V_j|=1$ for all $j\ne q$. 

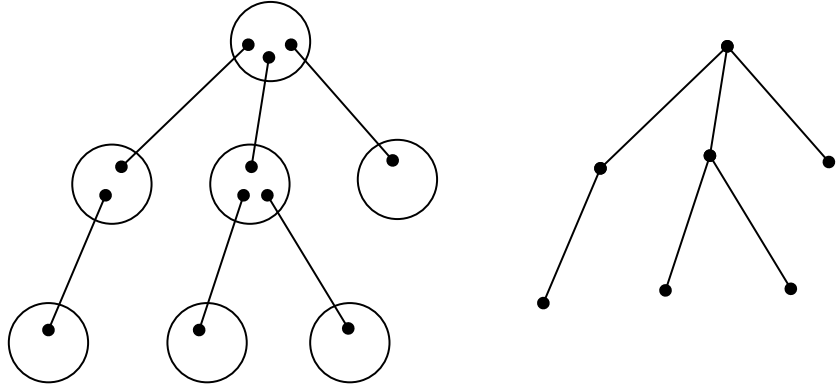
\begin{figure}[ht]
\centering

\tikzset{every picture/.style={line width=0.75pt,scale=0.8}} %set default line width to 0.75pt        

\begin{tikzpicture}[x=0.75pt,y=0.75pt,yscale=-1,xscale=1]
%uncomment if require: \path (0,303); %set diagram left start at 0, and has height of 303

%Shape: Circle [id:dp15980153815613596] 
\draw   (122,146) .. controls (122,132.19) and (133.19,121) .. (147,121) .. controls (160.81,121) and (172,132.19) .. (172,146) .. controls (172,159.81) and (160.81,171) .. (147,171) .. controls (133.19,171) and (122,159.81) .. (122,146) -- cycle ;
%Straight Lines [id:da9233244112955604] 
\draw    (233,58) -- (153,135) ;
\draw [shift={(153,135)}, rotate = 136.09] [color={rgb, 255:red, 0; green, 0; blue, 0 }  ][fill={rgb, 255:red, 0; green, 0; blue, 0 }  ][line width=0.75]      (0, 0) circle [x radius= 3.35, y radius= 3.35]   ;
\draw [shift={(233,58)}, rotate = 136.09] [color={rgb, 255:red, 0; green, 0; blue, 0 }  ][fill={rgb, 255:red, 0; green, 0; blue, 0 }  ][line width=0.75]      (0, 0) circle [x radius= 3.35, y radius= 3.35]   ;
%Shape: Circle [id:dp8568372079946195] 
\draw   (222,56) .. controls (222,42.19) and (233.19,31) .. (247,31) .. controls (260.81,31) and (272,42.19) .. (272,56) .. controls (272,69.81) and (260.81,81) .. (247,81) .. controls (233.19,81) and (222,69.81) .. (222,56) -- cycle ;
%Shape: Circle [id:dp7480788458282008] 
\draw   (209,146) .. controls (209,132.19) and (220.19,121) .. (234,121) .. controls (247.81,121) and (259,132.19) .. (259,146) .. controls (259,159.81) and (247.81,171) .. (234,171) .. controls (220.19,171) and (209,159.81) .. (209,146) -- cycle ;
%Shape: Circle [id:dp7102930778405164] 
\draw   (302,143) .. controls (302,129.19) and (313.19,118) .. (327,118) .. controls (340.81,118) and (352,129.19) .. (352,143) .. controls (352,156.81) and (340.81,168) .. (327,168) .. controls (313.19,168) and (302,156.81) .. (302,143) -- cycle ;
%Shape: Circle [id:dp17113949602220058] 
\draw   (82,246) .. controls (82,232.19) and (93.19,221) .. (107,221) .. controls (120.81,221) and (132,232.19) .. (132,246) .. controls (132,259.81) and (120.81,271) .. (107,271) .. controls (93.19,271) and (82,259.81) .. (82,246) -- cycle ;
%Shape: Circle [id:dp9057925268709781] 
\draw   (182,246) .. controls (182,232.19) and (193.19,221) .. (207,221) .. controls (220.81,221) and (232,232.19) .. (232,246) .. controls (232,259.81) and (220.81,271) .. (207,271) .. controls (193.19,271) and (182,259.81) .. (182,246) -- cycle ;
%Shape: Circle [id:dp3153426649913226] 
\draw   (272,246) .. controls (272,232.19) and (283.19,221) .. (297,221) .. controls (310.81,221) and (322,232.19) .. (322,246) .. controls (322,259.81) and (310.81,271) .. (297,271) .. controls (283.19,271) and (272,259.81) .. (272,246) -- cycle ;
%Straight Lines [id:da7388841313028028] 
\draw    (246,66) -- (235,135) ;
\draw [shift={(235,135)}, rotate = 99.06] [color={rgb, 255:red, 0; green, 0; blue, 0 }  ][fill={rgb, 255:red, 0; green, 0; blue, 0 }  ][line width=0.75]      (0, 0) circle [x radius= 3.35, y radius= 3.35]   ;
\draw [shift={(246,66)}, rotate = 99.06] [color={rgb, 255:red, 0; green, 0; blue, 0 }  ][fill={rgb, 255:red, 0; green, 0; blue, 0 }  ][line width=0.75]      (0, 0) circle [x radius= 3.35, y radius= 3.35]   ;
%Straight Lines [id:da07611902674521431] 
\draw    (260,58) -- (324,131) ;
\draw [shift={(324,131)}, rotate = 48.76] [color={rgb, 255:red, 0; green, 0; blue, 0 }  ][fill={rgb, 255:red, 0; green, 0; blue, 0 }  ][line width=0.75]      (0, 0) circle [x radius= 3.35, y radius= 3.35]   ;
\draw [shift={(260,58)}, rotate = 48.76] [color={rgb, 255:red, 0; green, 0; blue, 0 }  ][fill={rgb, 255:red, 0; green, 0; blue, 0 }  ][line width=0.75]      (0, 0) circle [x radius= 3.35, y radius= 3.35]   ;
%Straight Lines [id:da4476293239793685] 
\draw    (245,153) -- (296,237) ;
\draw [shift={(296,237)}, rotate = 58.74] [color={rgb, 255:red, 0; green, 0; blue, 0 }  ][fill={rgb, 255:red, 0; green, 0; blue, 0 }  ][line width=0.75]      (0, 0) circle [x radius= 3.35, y radius= 3.35]   ;
\draw [shift={(245,153)}, rotate = 58.74] [color={rgb, 255:red, 0; green, 0; blue, 0 }  ][fill={rgb, 255:red, 0; green, 0; blue, 0 }  ][line width=0.75]      (0, 0) circle [x radius= 3.35, y radius= 3.35]   ;
%Straight Lines [id:da6532006396155796] 
\draw    (230,153) -- (202,238) ;
\draw [shift={(202,238)}, rotate = 108.23] [color={rgb, 255:red, 0; green, 0; blue, 0 }  ][fill={rgb, 255:red, 0; green, 0; blue, 0 }  ][line width=0.75]      (0, 0) circle [x radius= 3.35, y radius= 3.35]   ;
\draw [shift={(230,153)}, rotate = 108.23] [color={rgb, 255:red, 0; green, 0; blue, 0 }  ][fill={rgb, 255:red, 0; green, 0; blue, 0 }  ][line width=0.75]      (0, 0) circle [x radius= 3.35, y radius= 3.35]   ;
%Straight Lines [id:da6658839634261211] 
\draw    (143,153) -- (107,238) ;
\draw [shift={(107,238)}, rotate = 112.95] [color={rgb, 255:red, 0; green, 0; blue, 0 }  ][fill={rgb, 255:red, 0; green, 0; blue, 0 }  ][line width=0.75]      (0, 0) circle [x radius= 3.35, y radius= 3.35]   ;
\draw [shift={(143,153)}, rotate = 112.95] [color={rgb, 255:red, 0; green, 0; blue, 0 }  ][fill={rgb, 255:red, 0; green, 0; blue, 0 }  ][line width=0.75]      (0, 0) circle [x radius= 3.35, y radius= 3.35]   ;
%Straight Lines [id:da7395053453528246] 
\draw    (535,59) -- (455,136) ;
\draw [shift={(455,136)}, rotate = 136.09] [color={rgb, 255:red, 0; green, 0; blue, 0 }  ][fill={rgb, 255:red, 0; green, 0; blue, 0 }  ][line width=0.75]      (0, 0) circle [x radius= 3.35, y radius= 3.35]   ;
\draw [shift={(535,59)}, rotate = 136.09] [color={rgb, 255:red, 0; green, 0; blue, 0 }  ][fill={rgb, 255:red, 0; green, 0; blue, 0 }  ][line width=0.75]      (0, 0) circle [x radius= 3.35, y radius= 3.35]   ;
%Straight Lines [id:da24132284051140684] 
\draw    (535,59) -- (524,128) ;
%Straight Lines [id:da9751122957671573] 
\draw    (535,59) -- (599,132) ;
\draw [shift={(599,132)}, rotate = 48.76] [color={rgb, 255:red, 0; green, 0; blue, 0 }  ][fill={rgb, 255:red, 0; green, 0; blue, 0 }  ][line width=0.75]      (0, 0) circle [x radius= 3.35, y radius= 3.35]   ;
\draw [shift={(535,59)}, rotate = 48.76] [color={rgb, 255:red, 0; green, 0; blue, 0 }  ][fill={rgb, 255:red, 0; green, 0; blue, 0 }  ][line width=0.75]      (0, 0) circle [x radius= 3.35, y radius= 3.35]   ;
%Straight Lines [id:da03817778782557946] 
\draw    (524,128) -- (575,212) ;
\draw [shift={(575,212)}, rotate = 58.74] [color={rgb, 255:red, 0; green, 0; blue, 0 }  ][fill={rgb, 255:red, 0; green, 0; blue, 0 }  ][line width=0.75]      (0, 0) circle [x radius= 3.35, y radius= 3.35]   ;
\draw [shift={(524,128)}, rotate = 58.74] [color={rgb, 255:red, 0; green, 0; blue, 0 }  ][fill={rgb, 255:red, 0; green, 0; blue, 0 }  ][line width=0.75]      (0, 0) circle [x radius= 3.35, y radius= 3.35]   ;
%Straight Lines [id:da9368509827664508] 
\draw    (524,128) -- (496,213) ;
\draw [shift={(496,213)}, rotate = 108.23] [color={rgb, 255:red, 0; green, 0; blue, 0 }  ][fill={rgb, 255:red, 0; green, 0; blue, 0 }  ][line width=0.75]      (0, 0) circle [x radius= 3.35, y radius= 3.35]   ;
\draw [shift={(524,128)}, rotate = 108.23] [color={rgb, 255:red, 0; green, 0; blue, 0 }  ][fill={rgb, 255:red, 0; green, 0; blue, 0 }  ][line width=0.75]      (0, 0) circle [x radius= 3.35, y radius= 3.35]   ;
%Straight Lines [id:da3949856450500677] 
\draw    (455,136) -- (419,221) ;
\draw [shift={(419,221)}, rotate = 112.95] [color={rgb, 255:red, 0; green, 0; blue, 0 }  ][fill={rgb, 255:red, 0; green, 0; blue, 0 }  ][line width=0.75]      (0, 0) circle [x radius= 3.35, y radius= 3.35]   ;
\draw [shift={(455,136)}, rotate = 112.95] [color={rgb, 255:red, 0; green, 0; blue, 0 }  ][fill={rgb, 255:red, 0; green, 0; blue, 0 }  ][line width=0.75]      (0, 0) circle [x radius= 3.35, y radius= 3.35]   ;

\end{tikzpicture}
\caption{An IMC in a 7-partite graph and the corresponding tree.}
\end{figure}

\begin{lm}[\textrm{\cite[Lemma~2.1]{HS2006}}] \label{it_fromimc}
	Let $G$ be an $r$-partite graph with vertex partition $V_1\cup \cdots \cup V_r$, and let $I=\{v_i,w_i:i\in [r-1]\}$ be an IMC in $G$, where $v_i$ and $w_i$ are adjacent. For any index $q\in [r]$, there is an $V_q$-avoiding PIT $T=\{s_j:j\in [r-1]\}$ in $G$ such that $s_j=v_j$ or $s_j=w_j$ for every $j\in [r-1]$.
\end{lm}

Given an IMC $I$ in $G$, define $A_v=A_v(I):=\{y\in V(G):N(y)\cap I=\{v\}\}$ for each $v\in I$. 

\begin{lm}[\textrm{\cite[Lemma~3.3]{HS2006}}]\label{rep_imc}
    Let $G\in \mathcal{G}_r(n)$ be a graph containing no IT and $I=\{a_i,b_i:i\in [r-1]\}$ be an IMC in $G$, where $a_i$ and $b_i$ are adjacent. Suppose vertices $v_i\in A_{a_i},w_i\in A_{b_i}$, $i\in [r-1]$ are independent in $G-E$, where $E=\{v_iw_i:i\in [r-1]\}$. Then $\{v_i,w_i:i\in [r-1]\}$ is an IMC in $G$.
\end{lm}

We call a graph $G\in \mathcal{G}_r(n)$ \textit{critical} if it contains no independent transversal but $G-e$ contains independent transversal for every edge $e\in E(G)$.

\begin{lm}[\textrm{\cite[Lemma~3.6]{HS2006}}]\label{edge_imc}
	If a graph $G\in \mathcal{G}_r(n)$ with $\Delta(G)<\frac{r-1}{2r-4}n$ ($r\ge 3$) is critical and $e\in E(G)$, then $e$ lies in an IMC in $G$.
\end{lm}

If $G\in \mathcal{G}_r(n)$ is critical with $\Delta(G) < \frac{r-1}{2r-4}n$ and $r\ge 7$, then \cite[Lemma~3.7]{HS2006} says that $G$ is the union of $r-1$ vertex-disjoint complete bipartite graphs. In fact, the proof of \cite[Lemma~3.7]{HS2006} is valid as long as $\Delta(G) < \frac{3r}{6r-7}n$, see the footnote below.\footnote{
\cite[Lemma~3.7]{HS2006} was proved with the help of \cite[Lemmas~3.4, 3.5]{HS2006}, and all three lemmas assume $\Delta = \Delta(G) < \frac{r-1}{2r-4}n$ and $r\ge 7$. In fact, each of the lemmas derives various lower bounds for $\Delta $, 
which contradict $\Delta < \frac{r-1}{2r-4}n$ when $r\ge 7$.
Below we list all these lower bounds and show that they all contradict $\Delta(G) < \frac{3r}{6r-7}n$ when $r\ge 2$.
    
%    \red{\cite[Lemma~3.7]{HS2006} is proved based on the assumption ``" and  which are also based on this assumption. We are going to explain why our assumption ``$\Delta(G) \le \frac{r}{2r-2}n+\varepsilon n$ with $\varepsilon \le \frac{1}{12r}$" gives contradictions as well in each place of these three lemmas the original assumption used. }
    
    In the proof of \cite[Lemma~3.4]{HS2006}, there are two lower bounds for $\Delta$. One is ``$\frac{|\mathscr{S}|2rn}{(4r-5)|\mathscr{S}|+1}\le \Delta$" for some set $\mathscr{S}$ with $|\mathscr{S}|\ge 3$, and the other is ``$2rn\le (4r-5)\Delta$".
    
    In the proof of \cite[Lemma~3.5]{HS2006}, there are four lower bounds for $\Delta$. The first is ``$(7-4r)\Delta +2rn \le 2\Delta$". The second is ``$|\mathscr{S}|2rn\le (4r-5)|\mathscr{S}|\Delta$". The third is ``$\frac{2r}{4r-5}n\le \Delta$". All these imply that $\Delta\ge \frac{2r}{4r-5}n$. The fourth is ``$\Delta\ge \frac{3r}{6r-7}n$".
    
    In the proof of \cite[Lemma~3.7]{HS2006}, there are three lower bounds for $\Delta$. The first and the third are ``$\Delta\ge \frac{3r}{6r-7}n$". The second is ``$\Delta\ge \frac{2r}{4r-5}n$".

    Note that $\frac{3r}{6r-7}n <\frac{2r}{4r-5}n$ for all $r\ge 2$. Thus, by assuming $\Delta < \frac{3r}{6r-7}n$, we obtain immediate contradictions with all the lower bounds except for 
     $\frac{|\mathscr{S}|2rn}{(4r-5)|\mathscr{S}|+1}\le \Delta$ for $|\mathscr{S}|\ge 3$. Indeed, since $\frac{|\mathscr{S}|2rn}{(4r-5)|\mathscr{S}|+1}$ is an increasing function of $|\mathscr{S}|$, we have $\frac{|\mathscr{S}|2rn}{(4r-5)|\mathscr{S}|+1}\ge \frac{3\cdot 2rn}{(4r-5)\cdot 3+1}=\frac{3rn}{6r-7}$, contradicting $\Delta < \frac{3r}{6r-7}n$ as well.
}

\begin{lm}[\textrm{\cite[Lemma~3.7]{HS2006}}]
\label{union_bipartite}
	If $G\in \mathcal{G}_r(n)$ is a critical graph with $\Delta(G) < \frac{3r}{6r-7}n$ ($r\ge 2$), then $G$ is the union of $r-1$ vertex-disjoint complete bipartite graphs.     %between $A_i$ and $B_i$, $i\in [r-1]$.
\end{lm}

When $G\in \mathcal{G}_r(n)$ is the union of $r-1$ vertex-disjoint complete bipartite graphs and $G$ is critical with $\Delta(G) < \frac{r-1}{2r-4}n$, we can apply Lemmas~\ref{rep_imc} and \ref{edge_imc} and show that $G$ contains many IMC's.

\begin{lm}\label{many_IMC}
    If $G\in \mathcal{G}_r(n)$ is a critical graph with $\Delta(G) < \frac{r-1}{2r-4}n$ ($r\ge 2$), and $G$ is the union of $r-1$ vertex-disjoint complete bipartite graphs between $A_i$ and $B_i$, $i\in [r-1]$, then $I=\{v_i,w_i:i\in [r-1]\}$ is an IMC in $G$ for any vertices $v_i\in A_i, w_i\in B_i$ where $i\in [r-1]$.
\end{lm}

\begin{proof}
    %We know that $G$ is the union of $r-1$ vertex-disjoint complete bipartite graphs between $A_i$ and $B_i$, $i\in [r-1]$.
    Fix an edge $e\in G$. Since $G$ is critical, by Lemma~\ref{edge_imc}, the edge $e$ lies in an IMC $I'=\{a_i,b_i:i\in [r-1]\}$, where $a_i$ and $b_i$ are adjacent. Since $G[I']$ is a matching of size $r-1$, for each $i\in [r-1]$, there is exactly one $j_i\in [r-1]$ such that $a_{j_i}$ and $b_{j_i}$ are in different parts of $G[A_i,B_i]$. After relabeling $a_1, b_1, \dots, a_{r-1}, b_{r-1}$ if necessary, we can assume $a_i\in A_i$ and $b_i\in B_i$ for all $i\in [r-1]$. Recall that $A_v=\{y\in V(G):N(y)\cap I'=\{v\}\}$ for $v\in I'$. We thus have $A_i=A_{b_i}$ and $B_i=A_{a_i}$. 
    Let $v_i\in A_i, w_i\in B_i$ for $i\in [r-1]$ and $I=\{v_i,w_i:i\in [r-1]\}$.
    Then $I$ is an independent set in $G-E$, where $E=\{v_iw_i:i\in [r-1]\}$.
    Then, by Lemma~\ref{rep_imc}, $I$ is an IMC in $G$. 
\end{proof}

%---------------------------------------------
\section{Proof of Theorem~\ref{mainresult} (lower bound)}
In this section we prove the lower bound in Theorem~\ref{mainresult}. Our proof indeed works for any $r\ge 3$ and any graph $G\in \mathcal{G}_r(n)$ with $\Delta(G)\le \frac{rn}{2r-2}- t$. 

Our proof starts with the following lemma, which shows that, if a graph $G\in \mathcal{G}_r(n)$ with $\Delta(G)\le \frac{rn}{2r-2}$ does not contain many independent transversals, then it is not far from the union of $r-1$ vertex-disjoint complete bipartite graphs with parts of size about $\frac{rn}{2r-2}$. We prove this lemma by applying the Graph Removal Lemma (Lemma~\ref{removal}) and the structure of critical graphs (Lemma~\ref{union_bipartite}).

	\begin{lm}\label{critical_graph}
		Given $r\ge 3$ and $0<\varepsilon \le 1/6$, there exists $\zeta=\zeta(r,\varepsilon)>0$ so that the following holds for sufficiently large $n$. %and all $1\le t\le rn/(2r-2)$. 
        Let $G\in \mathcal{G}_r(n)$ have parts $V_1,\cdots, V_r$ and satisfy $\Delta(G)\le \frac{r}{2r-2}n$. If $G$ contains at most $\zeta (rn)^r$ independent transversals, then there is an $r$-partite graph $G'=(V', E')$ with parts $V'_1\subseteq V_1, \dots, V'_r\subseteq V_r$ and the following properties.
        \begin{enumerate}
            \item $|V'_1| = \cdots = |V'_r|\ge n - (\eps n)/r$, 
            and $G'$ is a critical graph with $\Delta(G')\le \Delta(G)+\varepsilon n$. Furthermore, $\Delta(G'\setminus G)\le (\eps n)/r$ and $\Delta(G[V']\setminus G')\le (\eps n)/r$. 
            \item $G'$ is a union of $r-1$ vertex-disjoint complete bipartite graphs on $A_i\cup B_i$, $i\in [r-1]$ with 
            \begin{align}\label{eq:41}
                \frac{r n}{2r-2} - \eps n\le |A_i|, |B_i|\le \frac{r n}{2r-2} + \eps n.
            \end{align}
            \item$I=\{v_i,w_i:i\in [r-1]\}$ is an IMC in $G'$ for any vertices $v_i\in A_i, w_i\in B_i$ where $i\in [r-1]$.
        \end{enumerate}
  \end{lm}
	\begin{proof}    
    We first prove Part~1 and hold the proof of $\Delta(G[V']\setminus G')\le (\eps n)/r$ until the end because it requires Part~2.
    
		Let $\mu=\varepsilon^2/(16 r^6) $ and let $\zeta$ be $\zeta(r,\mu)$ from Lemma~\ref{removal}. Let $\tilde{G}=K(V_1,V_2,\cdots,V_r)\setminus G$.	Since $G$ contains at most $\zeta (rn)^r$ independent transversals, 
         the graph $\tilde{G}$ contains at most $\zeta(rn)^r$ copies of $K_r$. By Lemma~\ref{removal}, we can delete a set $E_1\subseteq E(\tilde{G})$ of at most $\mu (rn)^2$ edges from $\tilde{G}$ to make it $K_r$-free. 
        In other words, the graph $G_1 :=(V, E\cup E_1)$ is IT-free.
		
		Let $U_0=\{v\in V:d_{E_1}(v)\ge \frac{\varepsilon n}{4r^2}\}$. Since $|E_1|\le \mu (rn)^2$, it follows that 
        \[
        |U_0|\le 2\frac{\mu (rn)^2}{\varepsilon n/(4 r^2)}=\frac{\eps^2}{16 r^6} \frac{8 r^4n}{\varepsilon} =\frac{\varepsilon n}{2r^2}.
        \]
		Now we choose a set $U_1\supseteq U_0$ such that
        $|U_1\cap V_1| = \cdots =  |U_1\cap V_r| \le \varepsilon n/(2r^2)$. Let $V'=V(G)\setminus U_1$ and $V'_i = V_i\setminus U_1$ for $i\in [r]$.
        Then $n':=|V'_i| \ge n-\varepsilon n/(2r^2)$.
		Let $G_2=G_1[V']$. Then $G_2\in G_r(n')$ is IT-free with maximum degree 
		\[
        \Delta(G_2)\le \Delta(G)+\frac{\varepsilon n}{4r^2}\le \frac{rn}{2r-2}+\frac{\varepsilon n}{4 r^2}.
        \]
		
		Next we delete edges (one by one) from $G_2$ if necessary and obtain a spanning subgraph $G'$ of $G_2$ that is critical. Then $G'\in G_r(n')$ is
        an $r$-partite graph with parts $V'_1\dots, V'_r$    
        and $\Delta(G')\le \Delta(G_2)\le (\frac{r}{2r-2}+\frac{\varepsilon}{4r^2})n$. Since $V(G')=V'\subseteq V(G)\setminus U_0$, by the definition of $U_0$, we have 
        \begin{align}\label{eq:412}
            \Delta(G'\setminus G) \le \Delta(G_2\setminus G[V'])
        \le \frac{\eps n}{4r^2}.
        \end{align}
        This proves Part~1 except for $\Delta(G[V']\setminus G')\le \varepsilon n/r$. 
        In order to bound $\Delta(G[V']\setminus G')$, we need the structure of $G'$ given in Part~2 and
        thus work on Part~2 now.

        \medskip
        
	Since $\varepsilon \le 1/6$ and $r\ge 3$, it is easy to see that $\frac{r}{2r-2}+\frac{\varepsilon}{4r^2}\le \min \left\{\frac{3r}{6r-7},\frac{r-1}{2r-4}\right\}\left(1 -\frac{\eps}{2r^2}\right)$.\footnote{Indeed, $\frac{3r}{6r-7}-\frac{r}{2r-2}=\frac{r}{(6r-7)(2r-2)}>\frac{1}{12r}>\frac{\varepsilon}{4r^2}+\frac{3r}{6r-7}\frac{\varepsilon}{2r^2}$ and $\frac{r-1}{2r-4}-\frac{r}{2r-2}=\frac{1}{(r-2)(2r-2)}>\frac{1}{2r^2}>\frac{\varepsilon}{4r^2}+\frac{r-1}{2r-4}\frac{\varepsilon}{2r^2}$.} 
    This gives
    \[
    \Delta(G')\le \left(\frac{r}{2r-2}+\frac{\varepsilon}{4r^2}\right)n \le \left(\min \left\{\frac{3r}{6r-7},\frac{r-1}{2r-4}\right\}\right)\left(1 -\frac{\eps}{2r^2}\right) n \le \left(\min \left\{\frac{3r}{6r-7},\frac{r-1}{2r-4}\right\}\right)n'.
    \]   
		By Lemma~\ref{union_bipartite}, $G'$ is a union of $r-1$ vertex-disjoint complete bipartite graphs -- denote their parts by $A_i$ and $B_i$ for $i\in [r-1]$. By Lemma~\ref{many_IMC}, for any vertices $v_i\in A_i, w_i\in B_i$ where $i\in [r-1]$, the set $I=\{v_i,w_i:i\in [r-1]\}$ is an IMC in $G'$.
        
		Note that \[
        \max_{i\in[r-1]}\{|A_i|,|B_i|\}= \Delta(G')\le \left(\frac{r}{2r-2}+\frac{\varepsilon}{4r^2}\right)n.
        \]
        Since $v(G')=rn' \ge r (n - \frac{\eps n}{2r^2})$, we have
        
        \begin{align*}
			\min_{i\in[r-1]}\{|A_i|,|B_i|\}
            &\ge r n' -(2r-3)\max_{i\in[r-1]}\{|A_i|,|B_i|\}\\
			& \ge r \left(n - \frac{\eps n}{2r^2} \right) -(2r-3)\left(\frac{r n}{2r-2}+\frac{\varepsilon n}{4 r^2} \right) \\
			& = \frac{r n}{2r-2} - \frac{\eps n}{2r} - (2r-3)\frac{\eps n}{4r^2} = \frac{r n}{2r-2}-(4r-3)\frac{\varepsilon n}{4r^2}.
		\end{align*}
        Furthermore, by \eqref{eq:412},
        \begin{align*}
			\delta(G'\cap G) & \ge \delta(G')-\Delta(G'\setminus G) \ge \min_{i\in[r-1]}\{|A_i|,|B_i|\}-\frac{\eps n}{4r^2} \\
			& \ge \frac{r n}{2r-2}-(4r-3)\frac{\varepsilon n}{4r^2}-\frac{\eps n}{4r^2} \\
			& \ge \frac{r n}{2r-2}-\frac{\varepsilon n}{r}.
		\end{align*}

		We thus derive the desired bound
			\begin{align*}
            \hspace{2cm}
				\Delta(G[V']\setminus G') &\le \Delta(G)-\delta(G'\cap G) 
				\le \frac{r n}{2r-2} - \left(\frac{r n}{2r-2}-\frac{\varepsilon n}{r}\right)
				= \frac{\varepsilon n}{r}. \hspace{3cm} \qedhere
			\end{align*}
        \end{proof}	
        
	\medskip

\begin{proof}[Proof of Theorem \ref{mainresult} (lower bound)]	
     Let $r\ge 2$, $0< \varepsilon \ll 1$, $\zeta= \zeta(r, \eps)$ from Lemma~\ref{critical_graph}. Let $1\le t\le \zeta n$.
    Let $G\in \mathcal{G}_r(n)$ with maximum degree $\Delta(G)\le \lceil\frac{rn}{2r-2}\rceil-t$.
    Since 
    \[
    \zeta (rn)^r \ge \frac{r t}{2}\left(\frac{r n}{2r-2}\right)^{r-1} 
    \ge \frac{1}{4}\left( (t-1)(2r-2)+1 \right) \left(\frac{r n}{2r-2}\right)^{r-1}, 
    \]
    we may assume that $G$ contains at most $\zeta (rn)^r$ independent transversals (otherwise we are done). By Lemma~\ref{critical_graph}, there is an $r$-partite graph $G'$ with three properties stated in the lemma.

    Let $\tilde{E}= E(G[V']\setminus G') \cup E(G'\setminus G)$ consist of all edges of $G$ between $A_i\cup B_i$ and $A_j\cup B_j$ for all $i\ne j\in [r-1]$ and all non-edges of $G$ between $A_i$ and $B_i$ for all $i\in [r-1]$. By Lemma~\ref{critical_graph} Part~1, every vertex $v\in V'$ satisfies $d_{\tilde{E}}(v)\le \Delta(G[V']\setminus G')+\Delta(G'\setminus G)\le 2\varepsilon n/r$, and thus $|\tilde{E}|\le \frac{1}{2}|V'|(2\eps n/r)\le (r n)(\eps n/r)=\eps n^2$.

    For $i\in [r-1]$, by \eqref{eq:41}, we can find subsets $A'_i\subseteq A_i$ and $B'_i\subseteq B_i$ of size $m := \frac{r n}{2r-2} - \eps n$.
    By Lemma~\ref{critical_graph} Part~3,     
    \[
    \mathcal{I}:=\big\{\left\{v_1, w_1, \dots, v_{r-1}, w_{r-1}\right\}: v_i\in A_i', w_i\in B_i' \text{ for } i\in [r-1]\big\}   
    \]
    consists of $m^{2r-2}$ IMCs in $G'$, in which $v_i$ and $w_i$ are adjacent.
    Let $\I'\subseteq \I$ consist of all IMCs of $\I$ that induce no edge of $\tilde{E}$, in other words, every $I\in \I'$ satisfies $G[I]= G'[I]$. Since every edge of $G[A'_i\cup B'_i, A'_j\cup B'_j]$ for $i\ne j$ and every edge of $G'[A_i',B_i']\setminus G$ for $i\in [r-1]$ is contained in exactly $m^{2r-4}$ IMCs of $\I$, we have 
    \begin{align}\label{eq:42}
    |\I'| \ge |\mathcal{I}| - |\tilde{E}| m^{2r-4}\ge m^{2r-2}-\eps n^2 m^{2r-4}\ge \frac{3}{4}m^{2r-2}.
    \end{align}
    
    Arbitrarily take $I=\{v_i,w_i:i\in[r-1]\}\in \I'$.  
    Since $G[I]=G'[I]$ is a matching, we have 
    \begin{align*}
        |N_G(I)\cup I|&=|N_G(I)|\le \Delta(G)|I|\le \left(\Big \lceil\frac{rn}{2r-2}\Big \rceil-t\right)(2r-2) \\ 
        & \le \left(\frac{rn+2r-3}{2r-2}-t\right)(2r-2)= rn+(2r-3)-t(2r-2).
    \end{align*}    
    Consequently, $U_I:= V(G)\setminus (N_G(I)\cup I)$ has at least $t(2r-2)-(2r-3)=(t-1)(2r-2)+1$ vertices. Arbitrarily take a vertex $v\in U_I$ and suppose $v\in V_q$ for some $q\in [r]$. Since $I=\{v_i,w_i:i\in[r-1]\}$ is an IMC in $G'$, by Lemma~\ref{it_fromimc}, there is an $V_q$-avoiding PIT $T_q=\{s_j\in \{v_j,w_j\}:j\in [r-1]\}$ in $G'$. As $G[I]= G'[I]$, $T_q$ is also an $V_q$-avoiding PIT in $G$. Then $T=T_q\cup \{v\}$ is an IT of $G$. This process creates a triple $(I,v,T)$. Let $\mathcal{N}$ be the family of such triples and let $\mathcal{T}$ be the collection of (distinct) ITs in these triples. On one hand, by \eqref{eq:42}, we have 
    \begin{align}\label{eq:43}
    |\mathcal{N}|= \sum_{I\in \I'}\sum_{v\in U_I}1 \ge \frac{3}{4}m^{2r-2}\left((t-1)(2r-2)+1\right).
    \end{align}
    On the other hand, given $T\in \mathcal{T}$, we count the number of triples $(I,v,T)\in \mathcal{N}$ as follows. 
    By the process of producing $T$, we know that $|T\cap V'|\ge r-1$. If $|T\cap V'|=r-1$, then $v$ must be 
    the vertex in $T\setminus V'$. If $|T\cap V'|=r$, then 
    there exists a unique $i\in [r-1]$ such that $|T\cap (A'_i\cup B'_i)|=2$. Thus $v$ can be either of the two vertices in $T\cap (A'_i\cup B'_i)$. In either case we have at most two choices for $v$. Once $v$ is determined, there are at most $m^{r-1}$ choices to extend $T\setminus \{v\}$ to an IMC in $\I'$. Therefore, $|\mathcal{N}|\le |\mathcal{T}|\cdot 2 \cdot m^{r-1}.$
    Together with \eqref{eq:43}, this gives the desired bound
    \begin{align*}
        |\mathcal{T}| & \ge \frac{|\mathcal{N}|}{2m^{r-1}} \ge \frac{\frac{3}{4}m^{2r-2}\left((t-1)(2r-2)+1\right)}{2m^{r-1}} \\
        &= \frac{3}{8}\left((t-1)(2r-2)+1\right)m^{r-1} 
        \ge \frac{1}{4}\left((t-1)(2r-2)+1\right)\left(\frac{r n}{2r-2}\right)^{r-1},
    \end{align*}
    where the last inequality holds because $m^{r-1}=(\frac{rn}{2r-2}-\eps n)^{r-1}\ge \frac{2}{3}(\frac{rn}{2r-2})^{r-1}$ as $\eps\ll 1$.
 \end{proof}

%-------------------------------------------
\section{Proof of Theorem \ref{thm:blowup}}
In this proof, we will use a classical result of Erd\H os \cite{E1964} on the Tur\'an number of $K_r^r(s)$, \emph{the complete $r$-partite $r$-uniform hypergraph with $s$ vertices in each part}.

\begin{lm}[\cite{E1964}]\label{hyperclique}
	Given integers $r,s\ge 2$, let $n$ be sufficiently large. Then every $r$-uniform hypergraph on $n$ vertices with at least $n^{r-s^{1-r}}$ edges contains a copy of $K_r^r(s)$.
\end{lm}

In the proof of Lemma~\ref{hyperclique}, Erd\H{o}s applied the following result, which was derived from double counting and also appeared in \cite[Lemma~2.4]{BES1975}.

\begin{lm}[\cite{BES1975, E1964}]\label{intersectionsize}
	Let $n,k, q$ be natural numbers and $\frac{2q}{k}\le \alpha \le 1$. If $A_1,\cdots,A_k$ are subsets of $[n]$ such that $\frac{1}{k}\sum_{i=1}^k|A_i|\ge \alpha n$, then we can find distinct $A_{i_1},\cdots,A_{i_q}$ such that $|\bigcap_{j=1}^qA_{i_j}|\ge (\frac{\alpha}{2})^qn$.
\end{lm}

Recall that $z(n;s)$ is the maximum number of edges amongst the graphs $G\in \mathcal{G}_2(n)$ that do not contain $K_{s,s}$. K\H ov\'ari, S\'os and Tur\'an \cite{KST1954} gave the following upper bound for $z(n;s)$ (note that Lemma~\ref{intersectionsize} implies a slightly weaker bound). 

\begin{lm}[\cite{KST1954}]\label{findingKst}
	$z(n;s)\le (s-1)^{1/s}n^{2-1/s}+sn.$
\end{lm}

We are ready to prove Theorem \ref{thm:blowup}.

\begin{proof}[Proof of Theorem \ref{thm:blowup}]
    Let $C_{r,s}=4r^2s^{1/s}$, $t=C_{r,s}n^{1-1/s}$, and let $n$ be sufficiently large. Suppose $G\in \mathcal{G}_r(n)$ has parts $V_1, \dots, V_r$ and satisfies $\Delta(G)\le \frac{r}{2r-2}n - t$.
    
    If $r=2$, then the bipartite complement $\tilde{G}=K(V_1,V_2)\setminus G$ of $G$ satisfies $\delta(\tilde{G})\ge t$, and hence, $e(\tilde{G})\ge tn$. Since $t=4r^2s^{1/s}n^{1-1/s}$, we have $tn>(s-1)^{1/s}n^{2-1/s}+sn\ge z(n;s)$ by Lemma~\ref{findingKst}. Thus, $\tilde{G}$ contains a copy of $K_{s,s}$, which corresponds to an $IT(s)$ in $G$. 

    We assume $r\ge 3$ for the rest of the proof, where $r$ can be even or odd. 
	If $G$ contains at least $(rn)^{r-s^{1-r}}$ copies of IT, then we construct an $r$-uniform hypergraph $H$ on $V(G)$ whose hyperedges are all the ITs in $G$. Then $|E(H)|\ge (rn)^{r-s^{1-r}}.$ By Lemma~\ref{hyperclique}, $H$ contains a copy of $K_r^r(s)$, which corresponds an $IT(s)$ in $G$.

	We thus assume that $G$ contains less than $(rn)^{r-s^{1-r}}$ ITs. Let $\varepsilon \ll 1$, and $\zeta= \zeta(r, \eps)$ from Lemma~\ref{critical_graph}. Since $n$ is sufficiently large, we have $(rn)^{r-s^{1-r}}<\zeta (rn)^r$. By Lemma~\ref{critical_graph}, there is an $r$-partite graph $G'$ with three properties stated in the lemma.
	In particular, $G'=(V', E')$ consists of $r-1$ vertex-disjoint complete bipartite graphs $G'[A_i, B_i]$, $i\in [r-1]$ and $\Delta(G[V']\setminus G')\le \varepsilon n/r$.

    In our proof, we use the following claim. 
    \begin{clm}\label{find_weakIMCblowup}
       If there are sets $A_i'\subseteq A_i$ and $B_i'\subseteq B_i$ with $\min\{|A_i'|,|B_i'|\}\ge 2(i-1)s\varepsilon n/r+s$ for $i\in [r-1]$, then we can find $s$-element subsets $A_i''\subseteq A_i', B_i''\subseteq B_i'$, $i\in [r-1]$, such that for all $i\ne j\in [r-1]$, $G[A_i''\cup B_i'', A_j''\cup B_j'']=\emptyset$.    

    \end{clm}
    \begin{proof}
    We pick $A''_1, B''_1, \dots, A''_{r-1}, B''_{r-1}$ greedily as follows. Since $|A'_1|,|B'_1|\ge s$, we pick two arbitrary subsets $A''_1, B''_1$ of size $s$.
    Suppose we have chosen $A''_1, B''_1, \dots, A''_{i}, B''_{i}$ for some $i<r-1$. Note that any edge of $G$ between $A_i\cup B_i$ and $A_j\cup B_j$, $i\ne j$, is in $G[V']\setminus G'$.
    Since $\Delta(G[V']\setminus G')\le \varepsilon n/r$ and $A_1,B_1, \dots, A_{r-1}, B_{r-1}\subseteq V'$, we have 
    \[    
    \left|(A_{i+1}'\cup B_{i+1}') \cap N_G \left(A''_1\cup B''_1\cup \dots \cup A''_{i}\cup B''_{i}\right)\right|\le (2s i)\varepsilon n/r. 
    \]
    Since $|A'_{i+1}|,|B'_{i+1}|\ge (2i s)\varepsilon n/r+s$, we can pick two $s$-element subsets $A''_{i+1}, B''_{i+1}$ such that $G[A_j''\cup B_j'', A_{i+1}''\cup B_{i+1}'']=\emptyset$ for $j\le i$, as desired.
    \end{proof}

	Let $V_0=V(G)\setminus V'$ and $\alpha =1/(4s)$. Define 
    \[
    V_0'=\{v\in V_0:\overline{d}_G(v,A_i),\overline{d}_G(v,B_i)\ge \alpha n \text{ for every }i\in [r-1]\}
    \]
    as the set of vertices $v\in V_0$ with at least $\alpha n$ non-neighbors of $G$ in all $A_i$ and $B_i$ (including the vertices in the same part of $G$ as $v$).

\subsection*{Case 1: $|V_0'|\ge(\frac{2}{\alpha})^{2r-2}rs$.}   		

        We find the desired $IT(s)$ in three steps. Below we describe and justify each step.
    \begin{description}
        \item[Step 1.] Find subsets $S_0\subseteq V_0'\cap V_q$ of size $s$ for some $q\in [r]$,  $A'_i\subseteq A_i$ and $B'_i\subseteq B_i$ of size $|A'_i|, |B_i'|\ge \alpha' n$ such that $G[S_0, A'_i] = G[S_0, B'_i]= \emptyset$ for all $i\ge 1$, where $\alpha'= (\frac{\alpha}{2})^{(2/\alpha)^{2r-3}s}$.
    \end{description}

        Indeed, by averaging, $|V_0'\cap V_q|\ge (\frac{2}{\alpha})^{2r-2}s$  for some $q\in [r]$. Let $V_0^0$ be a subset of $V_0'\cap V_q$ of size $k:=(\frac{2}{\alpha})^{2r-2}s$. By the definition of $V_0'$, every vertex in $V_0^0$ has at least $\alpha n$ non-neighbors (of $G$) in $A_i$ and $B_i$ for all $i\in [r-1]$. Applying Lemma~\ref{intersectionsize} with $q=\frac{\alpha}{2}k$, we find a subset $V_0^1\subseteq V_0^0$ with $|V_0^1|=\frac{\alpha}{2}k$ such that the vertices in $V_0^1$ have at least $(\frac{\alpha}{2})^{|V_0^1|}n = \alpha' n$ common non-neighbors in $A_1$. Applying Lemma~\ref{intersectionsize} again, we find a subset $V_0^2\subseteq V_0^1$ with $|V_0^2|=\frac{\alpha}{2}|V_0^1|=(\frac{\alpha}{2})^2k$ such that vertices in $V_0^2$ have at least $(\frac{\alpha}{2})^{|V_0^2|}n\ge \alpha' n$ common non-neighbors in $B_1$ . Repeating this process, we obtain sets $S_0=V_0^{2r-2}\subseteq \cdots \subseteq V_0^1 \subseteq V_0^0$ such that $|S_0|=|V_0^{2r-2}|=(\frac{\alpha}{2})^{2r-2}k=s$ and the vertices in $S_0$ have at least $\alpha' n$ common non-neighbors in $A_i$ and $B_i$ for all $i\in [r-1]$. Denote these non-neighbors by $A_i'$ and $B_i'$ respectively. We have $|A'_i|, |B_i'|\ge \alpha' n$ and $G[S_0, A'_i] = G[S_0, B'_i]= \emptyset$. 

    \begin{description}
       \item[Step 2.] Find $s$-element subsets $A''_i\subseteq A'_i$ and $B''_i\subseteq B'_i$, $i\ge 1$ such that for all $i\ne j\in [r-1]$, $G[A''_i\cup B''_i, A''_j\cup B''_j]= \emptyset$, and for all $i\in [r]$, $A''_i\subseteq V_{a_i}$, $B''_i\subseteq V_{b_i}$ for some $a_i, b_i\in [r]$.
    \end{description}

        Indeed, since $|A_i'|,|B_i'|\ge \alpha' n \ge 2(i-1)s\varepsilon n+rs$ (by using $\varepsilon\ll 1$) for $i\in [r-1]$, by averaging, we can find subsets $\tilde{A}_i\subseteq A_i'\cap V_{a_i}, \tilde{B}_i\subseteq B_i'\cap V_{b_i}$ for some $a_i,b_i\in [r]$ with $|\tilde{A}_i|,|\tilde{B}_i|\ge 2(i-1)s\varepsilon n/r+s$. 
        Then, by Claim \ref{find_weakIMCblowup}, there are $s$-element sets $A_i''\subseteq \tilde{A}_i, B_i''\subseteq \tilde{B}_i, i\in [r-1]$ such that $G[A_i''\cup B_i'', A_j''\cup B_j'']=\emptyset$ for all $i\ne j$. 
        Note that $A_i''\subseteq V_{a_i},B_i''\subseteq V_{b_i}$ for $i\in [r-1]$. 

    \begin{description}
       \item[Step 3.] Choose $S_i= A''_i$ or $S_i= B''_i$ for $i\ge 1$ such that $S_0, S_1, \dots, S_{r-1}$ form an $IT(s)$.     
    \end{description}
			
		Indeed, we take a set $I=\{v_i\in A_i'',w_i\in B_i'':i\in [r-1]\}$.  By Lemma~\ref{critical_graph} Part~3, $I$ is an IMC in $G'$, in which $v_i$ and $w_i$ are adjacent for $i\in [r-1]$. 
        Then, by Lemma~\ref{it_fromimc}, $I$ contains an $V_q$-avoiding PIT $T=\{s_i:i\in [r-1]\}$ in $G'$, where $s_i=v_i$ or $s_i=w_i$ for $i\in [r-1]$. Let $S_i=A_i''$ if $s_i=v_i$ and $S_i=B_i''$ if $s_i=w_i$ for $i\in [r-1]$. For $i\ne j\in [r-1]$, since $G[A_i''\cup B_i'', A_j''\cup B_j'']=\emptyset$ from Step~2, it follows that $G[S_i, S_j]=\emptyset$. Together with $G[A_i'', S_0]=G[B_i'', S_0]=\emptyset$ and $S_0\subseteq V_q$ from Step~1, we conclude that $S_0,S_1,\cdots, S_{r-1}$ form an $IT(s)$ in $G$. 
        
 \medskip       
    This completes the proof of Case~1.

    \medskip
	In the rest of the proof, we assume $|V_0'|< (2/\alpha)^{2r-2}rs$. In particular, $|V_0'|<t=C_{r,s}n^{1-1/s}$. For $i\in [r-1]$, let \[
    U_i=\{v\in V_0:\overline{d}_G(v,A_i)< \alpha n\},\quad
    \text{and} \quad W_i=\{v\in V_0:\overline{d}_G(v,B_i)< \alpha n\}.
    \]
    Note that all $U_i$ and $W_i$ are disjoint because, say, $v\in U_1\cap W_1$ implies that 
    \begin{align*}
      d_G(v)&\ge d_G(v,A_1)+d_G(v,B_1)\ge \left(|A_1|-\overline{d}_G(v,A_1)\right)+\left(|B_1|-\overline{d}_G(v,B_1)\right)\\
    &\ge 2\left(\frac{r n}{2r-2}-\eps n-\alpha n\right) > \Delta(G) \quad \text{by \eqref{eq:41}},
    \end{align*}
    a contradiction.
    Thus, $U_1, W_1, \dots, U_{r-1}, W_{r-1}$ is a partition of $V_0\setminus V_0'$.

	Note that $A_1\cup W_1, B_1\cup U_1, \dots, A_{r-1}\cup W_{r-1}, B_{r-1}\cup U_{r-1}$ form a partition of $V(G)\setminus V'_0$. Without loss of generality, by averaging, assume $|B_1\cup U_1|\ge \frac{rn-|V_0'|}{2r-2} \ge \frac{rn-t}{2r-2}$. Since $\Delta(G)\le \frac{rn}{2r-2}-t$, each vertex $x\in A_1$ has at least $\frac{2r-3}{2r-2}t$ non-neighbors of $G$ in $B_1\cup U_1$. Thus,
    \[
    \overline{e}_G(A_1, B_1\cup U_1)\ge 
    \frac{2r-3}{2r-2}t\cdot |A_1| \ge \frac{2r-3}{2r-2}t\cdot \left(\frac{r}{2r-2}-\varepsilon\right)n\ge \frac{tn}{2}.
    \]
    It follows that $\overline{e}_G (A_1,B_1)\ge tn/4$ or $\overline{e}_G (A_1,U_1)\ge tn/4$. We proceed in these two cases.

\subsection*{Case 2: $\overline{e}_G (A_1,B_1)\ge tn/4$.}

    We find the desired $IT(s)$ in three steps that are similar to Case~1.
    \begin{description}
          \item[Step 1] Find subsets $A'_1\subseteq A_1$ and $B'_1\subseteq B_1$ of size $|A'_1|=|B'_1|=s$ such that $G[A'_1, B'_1] = \emptyset$ and $A'_1\in V_{a_1}$, $B'_1\subseteq V_{b_1}$ for some $a_1\ne b_1\in [r]$. 
    \end{description}
    
		Indeed, by averaging, there are $a_1,b_1\in [r]$ such that $\overline{e}_G (A_1\cap V_{a_1}, B_1\cap V_{b_1})\ge tn/(4r^2)$.
        Note that $a_1\ne b_1$. Indeed, since $G'[A_1, B_1]$ is complete bipartite by Lemma~\ref{critical_graph} Part 2, the two sets $A_1\cap V_{a_1}$ and $B_1\cap V_{b_1}$ are from different parts of $G'$. Thus, $V_{a_1}, V_{b_1}$ are two different parts of $G$ (Lemma~\ref{critical_graph} says that $G'$ follows the partition of $G$).     
        Recall that $|A_1|, |B_1|\le \frac{r n}{2r-2} + \eps n$ by \eqref{eq:41}. 
        Let $m:=\frac{r n}{2r-2} + \eps n$. Since $t=4r^2s^{1/s}n^{1-1/s}$ and $n$ is sufficiently large, we have  
        \[
        \overline{e}_G (A_1\cap V_{a_1}, B_1\cap V_{b_1})\ge tn/(4r^2) = s^{1/s}n^{2-1/s}>(s-1)^{1/s}m^{2-1/s}+sm\ge z(m;s)
        \]
        by Lemma~\ref{findingKst}. Thus, there are $s$-element subsets $A_1'\subseteq A_1\cap V_{a_1}$ and $B_1'\subseteq B_1\cap V_{b_1}$ such that $\overline{G}[A_1',B_1']$ is complete, which implies that $G[A_1', B_1']=\emptyset$.
        
    \begin{description}   
       \item[Step 2] Find $s$-element subsets $A'_i\subseteq A_i$ and $B'_i\subseteq B_i$, $i\ge 2$ such that for all $i\ne j\in [r-1]$, $G[A'_i\cup B'_i, A'_j\cup B'_j]= \emptyset$, and for all $i\ge 2$, $A'_i\subseteq V_{a_i}$, $B'_i\subseteq V_{b_i}$ for some $a_i, b_i\in [r]$.
    \end{description}

    The proof of this step is the same as the one in Case 1: for $i\ge 2$, by averaging, we find subsets $\tilde{A}_i\subseteq A_i\cap V_{a_i}, \tilde{B}_i\subseteq B_i\cap V_{b_i}$ for some $a_i,b_i\in [r]$ with $|\tilde{A}_i|,|\tilde{B}_i|\ge 2(i-1)s\varepsilon n/r+s$. 
    Then, by Claim \ref{find_weakIMCblowup}, we find $s$-element sets $A_i'\subseteq \tilde{A}_i, B_i'\subseteq \tilde{B}_i$ for $i\ge 2$ 
    (note that $A'_1$ and $B'_1$ were chosen in Step~1)
    such that $G[A_i'\cup B_i', A_j'\cup B_j']=\emptyset$ for all $i\ne j\in [r-1]$. 
  
    \begin{description}   
       \item[Step 3] Choose $S_i= A'_i$ or $S_i= B'_i$ for $i\ge 2$ such that $A'_1, B'_1, S_2, \dots, S_{r-1}$ form an $IT(s)$.     
    \end{description}
  
        Indeed, we take a set $I=\{v_i\in A_i',w_i\in B_i':i\in [r-1]\}$. By Lemma~\ref{critical_graph} Part~3, $I$ is an IMC in $G'$, in which $v_i$ and $w_i$ are adjacent for $i\in [r-1]$. Then, by Lemma~\ref{it_fromimc}, there is an $V_{a_1}$-avoiding PIT $T=\{s_i:i\in [r-1]\}$ in $G'$ with $s_i=v_i$ or $s_i=w_i$ for $i\in [r-1]$. Since $v_1\in A_1'\subseteq V_{a_1}$, we must have $s_1=w_1$. Let $S_i=A_i'$ if $s_i=v_i$ and $S_i=B_i'$ if $s_i=w_i$ for $i\in [r-1]$. In particular, $S_1=B_1'$. Since for all $i\ne j$, $G[A_i'\cup B_i', A_j'\cup B_j']=G[A_1', B_1']=\emptyset$ and $A_1'\subseteq V_{a_1}$, we have $A_1',B_1', S_2,\cdots, S_{r-1}$ form an $IT(s)$. This completes the proof of Case~2.
	
\subsection*{Case 3: $\overline{e}_G (A_1,U_1)\ge tn/4$.}	

    We find the desired $IT(s)$ in three steps, where Steps 2 and 3 are the same as in Case~1.
    \begin{description}
          \item[Step 1.] Find subsets $A'_1\subseteq A_1$ and $S_0\subseteq U_1$ of size $|A'_1|=|S_0|=s$ such that $G[A'_1, S_0] = \emptyset$ and $A'_1\subseteq V_{a_1}$, $S_0\subseteq V_q$ for some $a_1, q\in [r]$. 
    \end{description}
  
        Indeed, by averaging, there are $a_1,q\in [r]$ such that $\overline{e}_G (A_1\cap V_{a_1}, U_1\cap V_{q})\ge tn/(4r^2)$. Note that it is possible to have $a_1=q$.
		 Recall that $|A_1|\le (\frac{r}{2r-2}+\varepsilon)n$ and $|U_1|\le |V_0|\le \varepsilon n$ by Lemma~\ref{critical_graph} Part~2. 
         As shown in Case~2 Step~1, there are $s$-element subsets $A_1'\subseteq A_1\cap V_{a_1}$ and $S_0\subseteq U_1\cap V_{q}$ such that $\overline{G}[A_1',S_0]$ is complete, which implies that $G[A_1', S_0]=\emptyset$.

        Let $A'_i = A_i \setminus N_G(S_0)$, $i\ge 2$ and $B'_i = B_i \setminus N_G(S_0)$, $i\ge 1$. 
        Then $G[A'_i, S_0] = G[B'_i, S_0] = \emptyset$ for all $i\in [r-1]$.        
        The definition of $U_1$ implies that every $v\in U_1$ has at least $|A_1|- \alpha n$ neighbors in $A_1$ under $G$. Since $\Delta(G)\le \frac{r}{2r-2}n-t$, the vertex $v\in U_1$ has at most  
        \[
        \frac{r}{2r-2}n-t- \left(|A_1|- \alpha n\right)\le \frac{r}{2r-2}n-\left(\frac{r}{2r-2}n-\varepsilon n- \alpha n\right)\le (\varepsilon+\alpha)n
        \]
        neighbors in $V\setminus A_1$ under $G.$
        Therefore, $|A_i'|, 2\le i\le r-1$ and $|B_i'|,i\in [r-1]$ are at least 
        \[
        \min_{i\in [r-1]}\{|A_i|,|B_i|\}-s\cdot (\varepsilon +\alpha)n \ge \left(\frac{r}{2r-2}\right)n -\varepsilon n - s\cdot \left(\varepsilon +\frac{1}{4s}\right)n\ge \frac{n}{4}.
        \]

\medskip
       We omit Steps 2 and 3 because they are the same as in Case~1. This completes Case~3 and the proof of Theorem \ref{thm:blowup}.	
\end{proof}

%-------------------------------------
\section{Concluding Remarks}
\label{sec:cr}
In this paper we study $f_r(n, t)$, the minimum number of independent transversals in $r$-partite graphs with parts of size $n$ and maximum degree at most $\Delta_r(n) - t$. When $r$ is even and $t=o(n)$, our Theorem~\ref{mainresult} provides lower and upper bounds that differ by a factor of $4r^2$ when $t=1$ and a factor of $O(r)$ when $t\ge 2$. 

When $r$ is odd, Theorem~\ref{thm:threshold} says that $\Delta_r(n) = \lceil\frac{r-1}{2r-4}n\rceil$. Since $\frac{r}{2r-2}< \frac{r-1}{2r-4}$, it is more difficult to obtain an analogous result for odd $r$. In fact, Haxell and Szab\'o  \cite{HS2006} commented that both $f_r(n,1)=O(n^{r-2})$ and $f_r(n,1)=\Theta(1)$ were possible for odd $r$. Recall that it was shown \cite{BES1975} that $f_3(n, 1)=4$ for $n\ge 4$.

\medskip
Since $\frac12< \frac{r}{2r-2}< \frac{r-1}{2r-4}$, 
Theorems~\ref{thm:threshold} and \ref{thm:BES} together imply the following corollary.
Recall that a classical result of Haxell \cite{H2001} says that for any $r, n\in \mathbb{N}$, every $G\in G_r(n)$ with $\Delta(G)\le n/2$ contains an independent transversal.
\begin{col}\label{cor61}
    For any $r\ge 2$, there exist $c_r$ and $n_0$ such that every graph $G\in \mathcal{G}_r(n)$ with $n\ge n_0$ and 
	$\Delta(G)\le  n/2$ contains at least $c_r n^r$ independent transversals.
\end{col}

The aforementioned result of Wanless and Wood \cite{WW2022} has the following corollary.

\begin{col}[\textrm{\cite[Corollary 10]{WW2022}}]
\label{thm:ww}
    Let $r, n\in \mathbb{N}$.
   If $G\in G_r(n)$ has maximum degree $\Delta(G)\le n/4$, then $G$ contains at least $(n/2)^r$ independent transversals. 
\end{col}
Note that the degree condition in Corollary~\ref{thm:ww} is stronger than the one in Corollary~\ref{cor61} while Corollary~\ref{thm:ww} allows $r$ and $n$ to be arbitrary (for example, when $n$ is fixed and $r\to \infty$).
It is natural to ask whether $\Delta(G)\le n/4$ in Corollary~\ref{thm:ww} can be replaced by $\Delta(G)\le n/2$.

\medskip
We also obtain an asymptotically tight bound for the maximum degree of $r$-partite graphs containing no blowup of an independent transversal (the tightness assumes the widely believed lower bound for the Zarankiewicz number). 
In light of the aforementioned result \cite{MR5101285} for $r=3$, 
we believe that an analogous result holds for all odd $r$.
\begin{conj}\label{conj:TZ}
    For any odd number $r\ge 5$ and integer $s \ge 2$, there exist $C>0$ and $n_0\in \mathbb{N}$ such that every graph $G\in \mathcal{G}_r(n)$ with $n\ge n_0$ and $\Delta(G)\le \frac{r-1}{2r-4}n -Cn^{1-1/s}$ contains an $IT(s)$.
\end{conj}
By adding an isolated part of $n$ vertices to Construction~\ref{cons1}, the proof of Proposition~\ref{lem:lb2} shows that the degree condition in Conjecture~\ref{conj:TZ} cannot be weakened.

\medskip
Given integers $t\ge s\ge 1$, let $z(n;s,t)$ denote the largest number of edges in a bipartite graph of $\mathcal{G}_2(n)$ that contains no copy of $K_{s,t}$. By following the proof of Theorem~\ref{thm:blowup} closely and applying the well known bound $z(n;s,t)\le (t-1)^{1/s}n^{2-1/s}+sn$, we obtain a tight degree condition for $G\in \mathcal{G}_r(n)$ containing unbalanced blow-ups of independent transversal.

%Given an $r$-partite graph and integers $1\le s_1\le s_2\le \dots \le s_r$, an \textit{$(s_1,s_2,\dots,s_r)$-blowup of independent transversal}, denoted by $IT(s_1,s_2,\dots,s_r)$, is an independent set containing $s_1,s_2,\dots,s_r$ vertices from each part respectively. find $IT(s_1,s_r,\dots,s_r,s_r)$ when $\Delta(G)\le \frac{r}{2r-2}n- C_{r,s_1,s_r}n^{1-1/s_1}$. It implies a general version as follows.

\begin{thm}\label{thm:extra}
	For any even integer $r\ge 2$ and positive integers $s_1\le s_2\le \dots \le s_r$, there exist constants $C$ and $n_0$ such that every graph $G\in \mathcal{G}_r(n)$ with $n\ge n_0$ and maximum degree $\Delta(G)\le \frac{r}{2r-2}n- C n^{1-1/s_1}$ contains an $IT(s_1,\dots,s_r)$, an independent set with $s_1,\dots,s_r$ vertices in $r$ parts of $G$.  
\end{thm}
Below is a sketch of the proof of Theorem~\ref{thm:extra}. It suffices to show that, 
given integers $s_1\le s_r$, every graph $G\in \mathcal{G}_r(n)$ with sufficiently large $n$ and maximum degree $\Delta(G)\le \frac{r}{2r-2}n-4r^2s_r^{1/s_1}n^{1-1/s_1}$  contains an $IT(s_1,s_r,\cdots,s_r)$. 
The proof of this assertion is essentially the same as that of Theorem~\ref{thm:blowup}, except at a few points. First, if $G$ contains at least $(rn)^{r-s_r^{1-r}}$ ITs, then $G$ contains an $IT(s_r)$ by Lemma~\ref{hyperclique}. %on the hypergraph whose hyperedges are ITs in $G$. 
%If $G$ contains less than $(rn)^{r-s_r^{1-r}}$ ITs, then 
Otherwise we apply Lemma~\ref{critical_graph} and proceed to Case 1, where $|V'_0|\ge (\frac{2}{\alpha})^{2r-2}rs_1$ and $\alpha= 1/(4s_1)$. To handle this case, we follow the same three steps as in the proof of Theorem~\ref{thm:blowup}, except that in Step~2 we find $s_r$-element sets, and in Step~3 we find an $IT(s_1,s_r,\ldots,s_r)$. If $|V_0'|<(\frac{2}{\alpha})^{2r-2}rs_1$, then $t:=4r^2s_r^{1/s_1}n^{1-1/s_1}>(\frac{2}{\alpha})^{2r-2}rs_1>|V_0'|$ and we proceed to Cases~2 and~3 as in the proof of Theorem~\ref{thm:blowup}. The main difference occurs in Step~1 of these two cases, where we apply $z(m;s_1,s_r)\le (s_r-1)^{1/s_1}m^{2-1/s_1}+s_1 m$ and find a copy of $K_{s_1,s_r}$. 
%in $\overline{G}[A'_1, B'_1]$.
This is also why the exponent $1-1/s_1$ appears in Theorem~\ref{thm:extra}. 

%\textcolor{blue}{We can still do that and finish Case~1 by slightly changing the lower constant bound $(\frac{2}{\alpha})^{2r-2}rs$ of $|V_0'|$ to a new constant $(\frac{2}{\alpha})^{2r-2}rs_1$. In Step~2 we still applied Claim~\ref{find_weakIMCblowup} but to find $s_r$-elements sets, which is valid as $\varepsilon \ll 1$. If $|V_0'|<(\frac{2}{\alpha})^{2r-2}rs_1$, then we still have $t=4r^2s_r^{1/s_1}n^{1-1/s_1}>(\frac{2}{\alpha})^{2r-2}rs_1>|V_0'|$. Therefore, we can reduce the case to Case~2 and Case~3 as what we did in the proof of Theorem~\ref{thm:blowup}. The key step we discussed in these two cases is Step~1, where we applied the upper bound of $z(m;s)$ to find a $K_{s,s}$ in a dense bipartite graph. Now we can smoothly find a $K_{s_1,s_r}$ by applying the upper bound of $z(m;s_1,s_r)$ and given the higher density. It is also the reason why the exponent $1-1/s_1$ arises in the question. What left in these two cases are to find $s_r$-subsets in $A_i,B_i$ and close the argument of two cases as what we did in the proof of Theorem~\ref{thm:blowup}.}

The construction in Proposition~\ref{lem:lb2} shows that the term $n^{1-1/s_1}$ in Theorem~\ref{thm:extra} is tight if $z(n;s)= \Omega(n^{2-1/s})$. %the widely believed lower bound for the Zarankiewicz number holds.

\bigskip
\noindent\textbf{Acknowledgments.}
The authors thank the anonymous referees for valuable suggestions that improved the presentation of this paper.

\medskip
\noindent\textbf{Conflicts of interest.} 
None.

\medskip
\noindent\textbf{Financial support.}
Yantao Tang was partially supported by NSF grant DMS 2300346.
Yi Zhao was partially supported by NSF grant DMS 2300346 and Simons Collaboration Grant 710094.

\end{document}